\documentclass[12pt]{amsart}
\usepackage{amsmath}
\makeatletter
\newcommand{\leqnomode}{\tagsleft@true}
\newcommand{\reqnomode}{\tagsleft@false}
\makeatother
\usepackage{amsfonts}\usepackage{verbatim}
\usepackage{hyperref}
\usepackage{amssymb}\usepackage{marvosym}
\usepackage{dictsym}\usepackage{wasysym}
\usepackage{amstext}
\usepackage{amsbsy}\usepackage{bbold}
\usepackage{amsopn}\usepackage{MnSymbol}
\usepackage{amsthm}\usepackage{color}\usepackage{tikz}\usetikzlibrary{arrows}
\RequirePackage{filecontents}

\def\proclaim#1{\vskip0.5em\noindent{\bf #1}\it }
\def\endproclaim{\vskip0.5em\par\noindent\rm}

\def\proclaim#1{\vskip0.5em\noindent{\bf #1}\it}
\def\endproclaim{\vskip0.5em\par\noindent\rm}

\def\demo#1{\vskip0.5em\noindent{\bf #1\ }}

\def\text#1{\mbox{#1}}
\def\flushpar{\par\noindent}

\def\mod{\mbox{ mod }}

\newcommand{\mapright}[1]{%
    \smash{\mathop{%
        \hbox to 1cm{\rightarrowfill}
        }
    \limits^{#1}
    }
}
\newcommand{\mapleft}[1]{%
    \smash{\mathop{%
        \hbox to 1cm{\rightarrowfill}
        }
    \limits_{#1}
    }
}

\def\e{\epsilon}
\def\a{\alpha}
\def\b{\beta}
\def\G{\Gamma}
\def\g{\gamma}
\def\d{\delta}
\def\D{\Delta}
\def\s{\sigma}
\def\Si{\Sigma}
\def\th{\theta}
\def\l{\lambda}
\def\x{\times}

\def \F{\mathcal F}
\def\f{\flushpar}
\def\u{\underline}
\def\v{\varphi}

\def\om{\omega}

\def\B{\mathcal B}

\def\({\biggl(}
\def\){\biggr)}
\def\<{\langle}
\def\>{\rangle}

\def\bul{\smallskip\f$\bullet\ \ \ $}\def\lfl{\lfloor}\def\rfl{\rfloor}

\def\){\biggr)}

\def\<{\bold\langle}
\def\>{\bold\rangle}
\def\pprime{\prime\prime}

\def\bul{\smallskip\f$\bullet\ \ \ $}\def\sms{\smallskip\f}\def\sbul{\f$\bullet\
\ \ $}\def\Lra{\Longrightarrow}\def\Par{\smallskip\f\P}
\def\xbm{(X,\B,m)}

\begin{document}

\title{Discrepancy Skew Products and Affine  Random Walks}
\author{ Jon. Aaronson , Michael Bromberg $\text{and}$ Hitoshi Nakada}
 \address[Aaronson]{\ \ School of Math. Sciences, Tel Aviv University,
69978 Tel Aviv, Israel.\f\ \ \  {\it Webpage }: {\tt http://www.math.tau.ac.il/$\sim$aaro}}\email{aaro@post.tau.ac.il}
\address[Bromberg]{\ \ School of Math. Sciences, Tel Aviv University,r
69978 Tel Aviv, Israel.}\email{mic1@post.tau.ac.il}

\address[Nakada]{ Dept. Math., Keio University,Hiyoshi 3-14-1 Kohoku,
 Yokohama 223, Japan}\email[Nakada]{nakada@math.keio.ac.jp}
\thanks { The research of Aaronson and Bromberg was partially supported by ISF grant No. 1599/13. Nakada's research was partially supported by
JSPS grant No. 16K13766.\ \ \copyright  2016. }

\begin{abstract}We prove bounded rational ergodicity for some discrepancy skew products whose rotation number has bad rational approximation. This is done by considering the asymptotics of  associated affine random walks.

\end{abstract}
\maketitle\markboth{  
Discrepancy Skew Products and Affine  Random Walks.}{Aaronson, Bromberg,  $\text{and}$ Nakada}
\section*{\S1 Introduction}

\subsection*{Discrepancy skew products}
\par Let $\mathbb T:=\mathbb R/\mathbb Z\cong [0,1)$ denote the additive circle.

\par  
Consider  the  function $\v:\mathbb T\to\mathbb Z$ 
defined by
$$\v:=2\cdot 1_{[0,\tfrac12)}-1$$
and the {\it  skew products} 
$T_{\alpha}:\mathbb T\x\mathbb Z\to\mathbb T\x\mathbb Z$ 
defined for 
$\alpha\notin\mathbb Q$ by
$$T_\a(x,y)=(x+\a,y+\v(x)).$$ 
These are measure preserving transformations of the $\s$-finite measure space
$$\xbm=(\mathbb T\x\mathbb Z,\B(\mathbb T\x\mathbb Z),\text{\tt Leb}\x\#).$$

\

We have that  $T_\a^n(x,y)=(x+n\a,y+\v_n(x))$ where
$$\v_n(x)=\sum_{k=0}^{n-1}\v(x+k\a).$$
This is related to the discrepancy of the well-distribution of $(\{n\a\})_{n\ge 1}$ over $[0,\tfrac12)$ (see \cite{Schm}) and accordingly we call
the function $(n,x)\mapsto\v_n(x)$ the {\it discrepancy cocycle} and $T_\a$ the {\it discrepancy skew product}\ \ (aka
the
``cylinder flow arising from irregularity of distribution''  in \cite{Schm} $\text{and}$ ``deterministic random walk'' as in \cite{AK}). 

\

Ergodicity   of the discrepancy skew product $T_\a$ was established for $\a=\frac{\sqrt 5 -1}4$ in \cite{Schm} and then $\forall\ \a\in\mathbb R\setminus\mathbb Q$ in \cite{CK}.

\
\subsection*{Results}
\

We call $\a\in\mathbb R\setminus\mathbb Q$  {\it badly approximable} if
$$\inf\,\{q^2|\a-\tfrac{p}q|:\ q\in\mathbb N,\ p\in\mathbb Z\}>0$$
and denote
$$\text{\tt BAD}:=\{\a\in\mathbb R\setminus\mathbb Q:\ \a \ \ \text{badly approximable}\}.$$

Our main result is that if $\a\in\text{\tt BAD}$,
then $T_\a$ is {\it boundedly rationally ergodic} in the sense of \cite{BRE}. In particular, defining $\Psi_n=\Psi_n^{(\a)}:\mathbb T\to\mathbb N$ by
\begin{align*}
\Psi_n(x)& =S_n(1_{\mathbb T\x\{0\}})(x,0)\\ &:=\sum_{k=0}^{n-1}1_{\mathbb T\x\{0\}}\circ T_\a^k(x,0)\\ &=\#\{0\le k\le n-1:\ \v_k(x)=0\}\end{align*}
we prove that $\exists\ M>1$ so that
\begin{align*}\tag{\dsrailways}\label{choochoo}\|\Psi_n^{(\a)}\|_{L^\infty(\mathbb T)}\le M\int_\mathbb T \Psi_n^{(\a)}(t)dt\ \text{and}\ \ \int_\mathbb T \Psi_n^{(\a)}(t)dt=M^{\pm 1}\frac{n}{\sqrt{\log n}}. 
\end{align*}
\

Here and throughout, for $a,b>0,\ M>1$, $$a=M^{\pm 1}b\ \text{means}\ \ \frac1M\le \frac{a}b\le M.$$
\

We'll also consider more general subsequence versions of (\dsrailways):
\begin{align*}\tag{\dstechnical}\label{techy}\|\Psi_{\ell_n}^{(\a)}\|_{L^\infty(\mathbb T)}\le M\int_\mathbb T \Psi_{\ell_n}^{(\a)}(t)dt 
\end{align*}
which implies (as in \cite{BRE}) that there is a dense hereditary ring $R(T_\a)$ of sets of finite measure so that
$$\sum_{j=0}^{\ell_n-1}m(A\cap T^{-j}B)\ \sim\ m(A)m(B)a_{\ell_n}\  \text{as}\ n\to\infty\ \forall\ A,\ B\in R(T)$$ where
$a_{\ell_n}:=\int_\mathbb T \Psi_{\ell_n}^{(\a)}(t)dt$.

\

  In case $\a\in\mathbb R\setminus\mathbb Q$ quadratic, (\dsrailways)  was established in \cite{AK} and refined in \cite{ADDS}.
  \
  
  The discrepancy skew products occur as {\tt good sections} (in the sense of \cite{mode}) for directional translation flows of the infinite staircase translation surface
  as in  \cite{HHW} (see \cite{ADDS}).  Our result also holds for the corresponding directional, translation flow (with $\a$  badly approximable) 
  as can be seen using lemma 2.1 in \cite{mode}.
  \
  
\subsection*{Proof perspective}
\
 
The proof of (\dsrailways)  relies on a  weak, rough local limit theorem for an associated affine random walk arising from sequence of renormalizations related to the orbit of $2\a$ under a modified continued fraction transformation (see below). 
\

The quadratic case corresponds to a compact (eventually periodic) renormalization sequence, and the {\tt BAD} case 
corresponds to a precompact sequence.
  \subsection*{{\tt RAT}s $\text{and}$ {\tt ARW}s}
A {\it random affine transformation} ({\tt RAT}) on $\mathbb R^d$  is a random variable $F=(a(F),b(F))$ taking values in $M_{d\x d}(\mathbb R)\x\mathbb R^d$
(associated with the transformation
$x\mapsto ax+b$). 
\

We call the {\tt RAT} $F$ {\it discrete} if $F\in M_{d\x d}(\mathbb Z)\x\mathbb Z^d$ a.s.
\

An  {\it affine random walk} ({\tt ARW})  is a  $\mathbb R^d$-valued stochastic process $(X^{(n)})_{n\ge 0}$ defined by 
$$X^{(0)}:=0\ \text{and}\ X^{(n+1)}:=F_{n+1}(X^{(n)})=a{(F_n)}X^{(n)}+b{(F_n)}$$
where the $(F_n:\ n\ge 1)$  is a sequence of independent {\tt RAT}s  referred to as the {\it {\tt RAT} sequence}. 
In this paper, we only have need of a special kind of {\tt RAT} called {\tt flip type} (defined in \S4).
For other works on {\tt ARW}s (not of flip type), see \cite{FK},\ \cite{Kesten},  \cite{Vervaat}, \cite{GLP},  \cite{DI}  and references therein.
\

The rational ergodicity of the discrepancy skew product is governed by the asymptotic behavior of the temporal statistics of
the discrepancy cocycle values $(\v_n(0):\ n\ge 1)$ ( Visit Lemma 2.4 below). These temporal statistics are modeled by  certain {\tt ARW}s 
( Construction Lemma 4.1 below) and (\dsrailways) follows from
 a weak, rough,  local limit theorem  for the coordinates of these (Theorem 6.2 below).

  \section*{ \S2 Visit distributions of the discrepancy cocycle}
\

In this section, as in \cite{AK},  we show that (\dsrailways) (as on page \pageref{choochoo}) follows from certain asymptotic properties of ``{\tt visit distributions}'' (to be defined below).
\

We first note that we may assume without loss of generality that $0<\a<\frac12$.
\

This is because $-T_\a(-x,-n)=T_{1-\a}(x,n)$ and
$-\v_k^{(\a)}(x)=\v_k^{(1-\a)}(-x)$ whence

$$\Psi^{(\a)}_n(x)=\Psi^{(1-\a)}_n(-x).$$
Thus (\dsrailways) (as on page \pageref{choochoo}) for $\a\ \text{and}\ 1-\a$ are equivalent and we only consider the case $0<\a<\frac12$.
\subsection*{Calculation of the jump function orbit}
\

We recall from \cite{AK} the substitution algorithm to calculate the jump function orbit $(\v(\{n\a\}):\ n\ge 0)$.
\

We have, 
$$\v(n\a)=(-1)^{\sum_{j=1}^{n}\psi_j}$$
where
$$\psi_n=\psi_n^{(2\a)}:=1_{[1-2\a,1)}(2(n-1)\a).$$

   The {\it modified continued fraction expansion} of $\b\in (0,1)$ is
   $$\b=[n_1,n_2,\dots]=:1/n_1-1/n_2-1/n_3-\dots$$ with  each  $n_k\in\Bbb N_2:=\{a\in\Bbb N:\ a\ge 2\}$. Here (see \cite{KN}) $\b\in\Bbb Q$ iff $n_k\to 2$ as $n\to\infty$.   
   \proclaim{Theorem 2.1 \ \ (\cite{AK})}
   For $\b=[n_1,n_2,\dots]$, let
   $b_0(0)=0,\ b_0(1)=1\ \text{and}$
   $$b_{k+1}(0)=b_k(0)^{\odot (n_{k+1}-1)}\odot b_k(1)\ \ \text{and}\ \ b_{k+1}(1)=b_k(0))^{\odot (n_{k+1}-2)}\odot b_k(1),$$
   then
   $$(\psi^{(\b)}_1,\dots,\psi^{(\b)}_{\ell_k(0)})=b_k(0)\ \ \ (k\ge 1).$$\endproclaim
  Here $\odot$ denotes concatenation, and
   $\ell_k(i)=|b_k(i)|$ denotes the length of the block $b_k(i)\ \ (i=0,1)$.
   
   \
   
   For $\b=2\a=[n_1,n_2,\dots]$, set for $i=0,1$, $B_0(i):=[(-1)^i]$ and for $k\ge 0,\ i=0,1$ and $(n_{k+1},i)\ne (2,1)$:
     
  $$B_{k+1}(i)=\bigodot_{j=1}^{n_{k+1}-1-i}(-1)^{(j-1)\e_k(0)}B_k(0)\odot (-1)^{\e_k(0)(n_{k+1}-1-i)}B_k(1);$$
   where 
   $$\e_k(i):=\sum_{j=1}^{\ell_k(i)}(b_k(i))_j\ \mod 2$$ and
   $B_{k+1}(1)=B_k(1)$ in case $n_{k+1}=2$.
   \proclaim{Theorem 2.2  \ \ (\cite{AK})}
   $$B_k(0)=(\v(\{j\a\}))_{j=0}^{\ell_k(0)-1}\ \ \ (i=0,1).$$\endproclaim
  
\

\subsection*{Visit sets}
 \
 
 The {\it visit set to $\nu\in\mathbb Z$} is
 $$K_\nu:=\{n\ge 1:\ \v_n(0)=\nu\}$$
 and the {\it visit distributions} are the  measures $U_k^{(i)}$ on $\mathbb Z$ defined by
 $$U_k^{(i)}(\nu):=\#(K_\nu\cap [1,\ell_k(i)])\ \ \ (k\ge 1,\ i=0,1).$$
 \proclaim{Lemma 2.3}
 \begin{align*}
 &\tag{4.1}  \int_0^1\Psi_{\ell_k(0)}(x)dx\  \ge\ \ \frac1{4\ell_k(0)}\sum_{\nu\in\mathbb Z}[U_k^{(1)}(\nu)]^2;\\ &\tag{4.2}
 \int_0^1\Psi_{\ell_k(1)}(x)^Ndx\  \le\ \ \frac{2^N}{\ell_k(1)}\sum_{\nu\in\mathbb Z}[U_{k}^{(0)}(\nu)]^{N+1}\ \ \forall\ N\ge 1.
 \end{align*}
 \endproclaim

\

Statement (4.1) is  essentially Lemma 4.1 in \cite{AK} and proved in the same manner. Statement (4.2) is an upgrade of lemma 4.2 in \cite{AK}.
\demo{Proof of (4.2)}
\ \ As in the proofs of lemmas 4.1 and 4.2 in \cite{AK},
 
 \begin{align*}\ell_{k+r}(0)\int_0^1&\Psi_{\ell_k(1)}(x)^Ndx\\ &
 \underset{r\to\infty}{\text{\Large $\sim$}}
 \sum_{\nu\in\mathbb Z}\sum_{j\in [1,\ell_{k+r}(0)]\cap K_\nu}\#\,(K_\nu\cap [j+1,j+\ell_k(1)])^N\end{align*}
and for $r\ge 1,\ \exists\ J=J_{r,k}\ge 1,\ 1=m_1<\dots<m_J$ and $\e_1,\dots\e_{J-1}=\pm 1,\ i_1,\dots,i_{J-1}=0,1$  so that 
$$m_{j+1}-m_j=\ell_k(i_j)\ \forall\ j,\ \ [1,\ell_{k+r}(0)]=\bigcupdot_{j=1}^{J-1}[m_j,m_{j+1}),$$ 
$$(\v(m_j\a),\v(m_j+1)\a),\dots,\v((m_{j+1}-1)\a))=\e_jB_k(i_j).$$

\

For fixed $\nu\in\mathbb Z,$
 \begin{align*}
  K_\nu\cap [m_j,m_{j+1})=m_j+K_{\e_j(\nu-\v_{m_j}(0))}\cap [1,\ell_k(i_j)].
 \end{align*}

Note that 
 $$(\v(m_j\a),\v((m_j+1)\a),\dots,\v((m_{j+1}+\ell_k(1)-1)\a))=(\e_jB_k(i_j),\D_jB_k(1))$$
 for some $\D_j=\pm 1$.

 We have 
as before, for fixed $\nu\in\mathbb Z$,
\begin{align*}
\sum_{i\in [1,\ell_{k+r}(0)]\cap K_\nu}\#\,(K_\nu\cap [i+1,i+\ell_k(1)])^N&=
\sum_{j=1}^{J-1}\sum_{i\in [m_j,m_{j+1})\cap K_\nu}\#\,(K_\nu\cap [i+1,i+\ell_k(1)])^N\\ &\le
\sum_{j=1}^{J-1}\sum_{i\in [m_j,m_{j+1})\cap K_\nu}\#\,(K_\nu\cap [i+1,m_{j+1}+\ell_k(1))])^N.
\end{align*}
Fix $j$. For fixed $i\in [m_j,m_{j+1})$,
\begin{align*}\#\,(K_\nu\cap [i+1,m_{j+1}+\ell_k(1))])&=\#\,(K_\nu\cap [i+1,m_{j+1}])+\#\,(K_\nu\cap [m_{j+1}+1,m_{j+1}+\ell_k(1))])\\ &=
\#\,(K_\nu\cap [i+1,m_{j+1}])+U_k^{(1)}(\D_j(\nu-\v_{m_{j+1}(0)})). 
\end{align*}
Thus 
\begin{align*}
 &\sum_{i\in [m_j,m_{j+1})\cap K_\nu}\#\,(K_\nu\cap [i+1,m_{j+1}+\ell_k(1))])^N=
 \\ &=\sum_{i\in [m_j,m_{j+1})\cap K_\nu}(\#\,(K_\nu\cap [i+1,m_{j+1}])+U_k^{(1)}(\D_j(\nu-\v_{m_{j+1}(0)})))^N
 \\ &=\sum_{r=0}^N\binom{N}r\(\sum_{i\in [m_j,m_{j+1})\cap K_\nu}\#\,(K_\nu\cap [i+1,m_{j+1}])^r\)U_k^{(1)}(\D_j(\nu-\v_{m_{j+1}(0)}))^{N-r}
 \\ &\le\sum_{r=0}^N\binom{N}r U_k^{(i_j)}(\e_j(\nu-\v_{m_{j}(0)}))^{r+1}U_k^{(1)}(\D_j(\nu-\v_{m_{j+1}(0)}))^{N-r}
 \\ &\le\sum_{r=0}^N\binom{N}r U_k^{(0)}(\e_j(\nu-\v_{m_{j}(0)}))^{r+1}U_k^{(0)}(\D_j(\nu-\v_{m_{j+1}(0)}))^{N-r}.
\end{align*}\

Using this and H\"older's inequality,
\begin{align*}
 &\sum_{\nu\in\mathbb Z}\sum_{j\in [1,\ell_{k+r}(0)]\cap K_\nu}\#\,(K_\nu\cap [j+1,j+\ell_k(1)])^N\le 
 \\ &\le\sum_{j=1}^J\sum_{r=0}^N\binom{N}r\sum_{\nu\in\mathbb Z} U_k^{(0)}(\e_j(\nu-\v_{m_{j}(0)}))^{r+1}U_k^{(0)}(\D_j(\nu-\v_{m_{j+1}(0)}))^{N-r}
  \\ &\le\sum_{j=1}^J\sum_{r=0}^N\binom{N}r\(\sum_{\nu\in\mathbb Z} U_k^{(0)}(\e_j(\nu-\v_{m_{j}(0)}))^{N+1}\)^{\frac{r+1}{N+1}}
  \(\sum_{\nu\in\mathbb Z}U_k^{(0)}(\D_j(\nu-\v_{m_{j+1}(0)}))^{N+1}\)^{\frac{N-r}{N+1}}
  \\ &\le J\sum_{r=0}^N\binom{N}r\(\sum_{\nu\in\mathbb Z} U_k^{(0)}(\nu)^{N+1}\)^{\frac{r+1}{N+1}}
  \(\sum_{\nu\in\mathbb Z}U_k^{(0)}(\nu)^{N+1}\)^{\frac{N-r}{N+1}}
  \\ &=2^NJ\sum_{\nu\in\mathbb Z}U_k^{(0)}(\nu)^{N+1}
  \end{align*}
whence
\begin{align*}\int_0^1\Psi_{\ell_k(1)}(x)^Ndx&
 \xleftarrow[r\to\infty]{}\frac1{\ell_{k+r}(0)}\sum_{\nu\in\mathbb Z}\sum_{j\in [1,\ell_{k+r}(0)]\cap K_\nu}\#\,(K_\nu\cap [j+1,j+\ell_k(1)])^N
 \\ &\le \frac{2^NJ}{\ell_{k+r}(0)}\sum_{\nu\in\mathbb Z}U_k^{(0)}(\nu)^{N+1}
 \\ &\le\frac{2^N}{\ell_{k}(1)}\sum_{\nu\in\mathbb Z}U_k^{(0)}(\nu)^{N+1}.\ \ \ \CheckedBox\text{\rm (4.2)}
 \end{align*}

\sms{\bf Visit lemma}
\

This is a Fourier series (or generating function) consequence of   Lemma 2.3.   
 
 \
 
 Let
 $$\widehat{U}_k^{(i)}(Z):=\sum_{\nu\in\mathbb Z}U_k^{(i)}(\nu)Z^\nu\ \ \ \  (i=0,1,\ |Z|=1).$$
 
\proclaim{ Visit lemma 2.4} 
\begin{align*}&\tag{4.1'} \int_0^1\Psi_{\ell_k(0)}(x)dx\  \ge\ \ \frac1{4\ell_k(0)}\int_{\mathbb S^1}|\widehat{U}_k^{(1)}(Z)|^2dZ
\\ &\tag{4.2'}\|\Psi_{\ell_k(1)}\|_\infty\le 2\int_{\mathbb S^1}|\widehat{U}_k^{(1)}(Z)|dZ.
\end{align*}\endproclaim
Where here and throughout,
$$\int_{\mathbb S^1}f(Z)dZ:=\frac1{2\pi}\int_0^{2\pi}f(e^{i\th})d\th.$$
\demo{Proof}\ \ For fixed $k\ge 1$,
\

The statement (4.1') follows from (4.1) via the Riesz-Fischer theorem and the statement 
 (4.2'),  follows from (4.2) using the Hausdorff-Young theorem as in the proof of theorem 6.1 of \cite{AK}.\ \ \Checkedbox

\section*{\S3  Visit distribution  transitions ({\rm as in \cite{AK}})}
       Set  $s_k(i):=\sum_{j=0}^{\ell_k(i)-1}\v(\{j\a\})$ and define the {\it orbit blocks}
 $$\Si_k(i):=(\v_1(0),\v_2(0),\dots,\v_{\ell_k^{(i)}}(0)).$$
Our goal here is to obtain  the  transitions of the visit distribution generating functions.
 \sms{\bf Transitions in terms of blocks}
 \
 
 From theorem 2.2 above we see that for $k\ge 1,\ i=0,1$, where $(n_{k+1},i)\ne (2,1)$:
 \bul in case $\e_k(0)=0$,
 $$\Si_{k+1}(i)=\bigodot_{j=1}^{n_{k+1}-1-i}(\Si_k(0)+(j-1)s_k(0)\mathbb{1})\odot(\Si_k(1)+(n_{k+1}-i-1)s_k(0)\mathbb{1});$$
 \bul in case $\e_k(0)=1\ \text{and}\ \ n_{k+1}-1-i\in 2\mathbb Z$,
  $$\Si_{k+1}(i)=[\Si_k(0),s_k(0)\mathbb{1}-\Si_k(0)]^{\odot\frac{n_{k+1}-1-i}2}\odot\Si_k(1)$$
  \bul in case $\e_k(0)=1\ \text{and}\ \ n_{k+1}-1-i\in 2\mathbb Z+1$,
  $$\Si_{k+1}(i)=[\Si_k(0),s_k(0)\mathbb{1}-\Si_k(0)]^{\odot\frac{n_{k+1}-2-i}2}\odot\Si_k(0)\odot (s_k(0)\mathbb{1}-\Si_k(1)).$$
 \bul  and
   $\Si_{k+1}(1)=\Si_k(1)$ in case $n_{k+1}=2$.
   \

  \sms{\bf Transitions in terms of  visit distributions}
  \
  
   For $k\ge 1,\ i=0,1\ \text{and}\ \nu\in\mathbb Z$, where $(n_{k+1},i)\ne (2,1)$:
  \bul in case $\e_k(0)=0$,
 $$U_{k+1}^{(i)}(\nu)=\sum_{j=1}^{n_{k+1}-1-i}U_{k}^{(0)}(\nu-(j-1)s_k(0))+U_{k}^{(1)}(\nu-(n_{k+1}-i-1)s_k(0));$$
 \bul in case $\e_k(0)=1\ \text{and}\ \ n_{k+1}-1-i\in 2\mathbb Z$,
  $$U_{k+1}^{(i)}(\nu)=\frac{n_{k+1}-1-i}2 (U_{k}^{(0)}(\nu) + U_{k}^{(0)}(s_k(0)-\nu))+U_{k}^{(1)}(\nu);$$
  \bul in case $\e_k(0)=1\ \text{and}\ \ n_{k+1}-1-i\in 2\mathbb Z+1$,
  $$U_{k+1}^{(i)}(\nu)=\frac{n_{k+1}-i}2 U_{k}^{(0)}(\nu) + \frac{n_{k+1}-2-i}2U_{k}^{(0)}(s_k(0)-\nu))+U_{k}^{(1)}(s_k(0)-\nu);$$
 \bul  and
   $U_{k+1}^{(1)}(\nu)=U_{k}^{(1)}(\nu)$ in case $n_{k+1}=2$.

   \

   \sms{\bf Transitions in terms of  generating functions}
   \
   
   For $k\ge 1,\ i=0,1\ \text{and}\ Z\in\mathbb S^1$, where $(n_{k+1},i)\ne (2,1)$:
  \bul in case $\e_k(0)=0$,
 $$\widehat{U}_{k+1}^{(i)}(Z)=\sum_{j=1}^{n_{k+1}-1-i}Z^{(j-1)s_k(0)}\widehat{U}_{k}^{(0)}(Z)+Z^{(n_{k+1}-i-1)s_k(0)}\widehat{U}_{k}^{(1)}(Z);$$
 \bul in case $\e_k(0)=1\ \text{and}\ \ n_{k+1}-1-i\in 2\mathbb Z$,
  $$\widehat{U}_{k+1}^{(i)}(Z)=\frac{n_{k+1}-1-i}2 (\widehat{U}_{k}^{(0)}(Z) + Z^{s_k(0)}\widehat{U}_{k}^{(0)}(Z^{-1})+\widehat{U}_{k}^{(1)}(Z);$$
  \bul in case $\e_k(0)=1\ \text{and}\ \ n_{k+1}-1-i\in 2\mathbb Z+1$,
  $$\widehat{U}_{k+1}^{(i)}(Z)=\frac{n_{k+1}-i}2 \widehat{U}_{k}^{(0)}(Z) + \frac{n_{k+1}-2-i}2Z^{s_k(0)}\widehat{U}_{k}^{(0)}(Z^{-1})+Z^{s_k(0)}\widehat{U}_{k}^{(1)}(Z^{-1});$$
 \bul  and
   $\widehat{U}_{k+1}^{(1)}(Z)=\widehat{U}_{k}^{(1)}(Z)$ in case $n_{k+1}=2$.

\subsection*{Simplified visit distributions}
\

 Our next task is to eliminate the dependence of the visit distribution transitions on the positions $s_k(0)$. We'll do this exploiting the
parity sequence $(\u\e_k:\ k\ge 1)$.
\sms{\bf Parities, block lengths and positions}
\

 Let $k\ge 1$. The $k^{\text{\tiny th}}$:
 \begin{align*}  
&\text{\it parity state vector  {\rm is}}\ \  \u\e_k:=\left(\begin{matrix}  \e_k(0) 
\\  \e_k(1)\end{matrix}\right);\\ &\text{\it block length vector {\rm is}} \ \ \u\ell_k:=\left(\begin{matrix}  \ell_k(0) 
\\  \ell_k(1)\end{matrix}\right),\\ &\text{\it position vector {\rm is}} \ \ \ \ \u s_k:=\left(\begin{matrix}  s_k(0) 
\\  s_k(1)\end{matrix}\right) \end{align*}
where $s_k(i):=\sum_{j=0}^{\ell_k(i)-1}\v(\{j\a\})$.

Next let 
 \begin{align*} &A(0):=\left(\begin{matrix}  1  & -1 
\\  0  & 1  \end{matrix}\right)\ \ A(1):=\left(\begin{matrix}  0  & 1 
\\  1  & -1  \end{matrix}\right)\\ & B(n):=\left(\begin{matrix}  n-1  & 1 
\\  n-2  & 1  \end{matrix}\right)\ \ \text{and}\ A(n):=A(n\ \mod 1),  
 \end{align*}

then
$$\u s_{k+1}=\begin{cases}& B(n_{k+1})\u s_k\ \ \ \ \ \ \e_k(0)=0;\\
& A(n_{k+1})\u s_k\ \ \ \ \ \ \e_k(0)=1.              
             \end{cases}$$
             
             The following diagram gives (as in \cite{AK}) the parity transitions $\u\e_k\to\u\e_{k+1}$ conditional on $N=n_{k+1}$ together with the corresponding position
             vector transitions.
             
             \

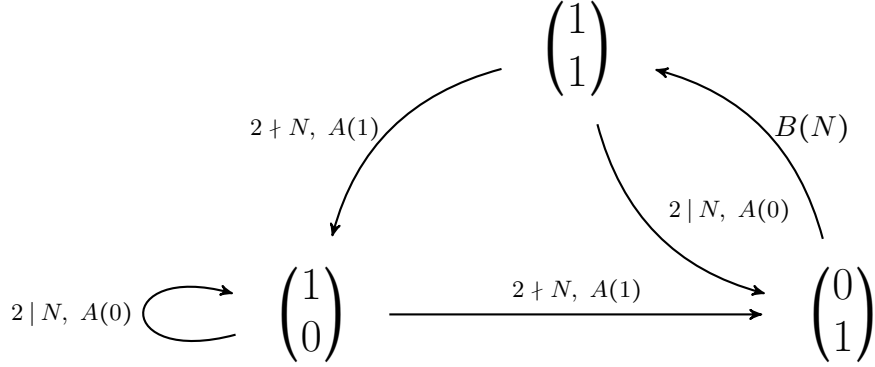
\begin{tikzpicture}[->,>=stealth',shorten >=1pt,auto,node distance=5cm,
  thick,main node/.style={circle,draw=none,font=\sffamily\Large\bfseries}]
  \node[main node] (1) {$\left(\begin{matrix}  1 \\  1\end{matrix}\right)$};
  \node[main node] (2) [below left of=1] {$\left(\begin{matrix}  1 \\  0\end{matrix}\right)$};
  \node[main node] (4) [below right of=1] {$\left(\begin{matrix}  0 \\  1\end{matrix}\right)$};
  \path[every node/.style={font=\sffamily\small}]
    (1) edge [bend right] node [bend left]  {\scriptsize $2\mid N,\  A(0)$} (4)
        edge [bend right] node[left] {\scriptsize$2\nmid N,\ A(1)$} (2)
   (2)  edge  node  {\scriptsize$2\nmid N,\ A(1)$} (4)
          edge [loop left] node {\scriptsize $2\mid N,\  A(0)$} (2)
    (4)  edge [bend right] node[right] {$B(N)$} (1);
\end{tikzpicture}
 \proclaim{Lemma 3.1 (\cite{AK})}\ \ \ For $k\ge 1$,
 $$s_k(0)=1,\ \ s_k(1)=1\ \ \text{or}\ \ s_k(0)-s_k(1)=1$$
 according to whether
 $$\u\e_k:=\left(\begin{matrix}  0 \\  1\end{matrix}\right),\ \left(\begin{matrix}  1 \\  0\end{matrix}\right),\ \left(\begin{matrix}  1 \\  1\end{matrix}\right)\ \ \text{respectively}.$$
 \endproclaim

\

To this end, define $T_k=T_k(\u\e_k)=T_k(\e_k(0))$ by:

\begin{align*}T_k\left(\begin{matrix}  1\\  0\end{matrix}\right)= T_k\left(\begin{matrix}  1\\  1\end{matrix}\right)=-s_k(0),\ \ \ 
T_k\left(\begin{matrix} 0\\  1\end{matrix}\right)=-s_k(1).
\end{align*}
Evidently 
\begin{align*}
\tag{\dsmilitary}\u\e_k(0)=1\ \overset{\text{\tiny{$T_k=-s_k(0)$}}}\Lra\ T_{k+1}- T_k= T_{k+1}+ T_k+2s_k(0)
\end{align*}
The following is established by straightrforward computation.
\proclaim{Proposition 3.2: \ {\rm Increments of the $T_k$s}}
\begin{align*}&\tag{1}T_{k+1}\left(\begin{matrix}  0\\  1\end{matrix}\right)- T_k\left(\begin{matrix}  1\\  0\end{matrix}\right)=1,
\\  &\tag{2}T_{k+1}\left(\begin{matrix}  1\\  0\end{matrix}\right)- T_k\left(\begin{matrix}  1\\  0\end{matrix}\right)=
1,
\\  &\tag{3}T_{k+1}\left(\begin{matrix}  1\\  0\end{matrix}\right)- T_k\left(\begin{matrix}  1\\  1\end{matrix}\right)=1
\\  &\tag{4}T_{k+1}\left(\begin{matrix}  0\\  1\end{matrix}\right)- T_k\left(\begin{matrix}  1\\  1\end{matrix}\right)=1
\\  &\tag{5}T_{k+1}\left(\begin{matrix}  1\\  1\end{matrix}\right)- T_k\left(\begin{matrix}  0\\  1\end{matrix}\right)=
-(n_{k+1}-1).\end{align*}\endproclaim

\subsection*{Simplified visit distribution transitions}
\

Given $\a\in\mathbb T\setminus\mathbb Q$ we define the {\it simplified visit distributions} by
$$V_k^{(i)}(J):=\begin{cases}&U_k^{(i)}(\nu)\ \ \ \ \ \ J=2\nu+T_k,\ \nu\in\mathbb Z,\\ &0\ \ \ \ \ \ \text{else.}\end{cases}$$
where $T_k=T_k(\e_k(0))$.
\

The generating functions are given by
$$\widehat{V}_k^{(i)}(Z)=Z^{T_k}\widehat{U}_k^{(i)}(Z^2).$$
\

They have the property that
$$\int_{\mathbb S^1}|\widehat{V}_k^{(i)}(Z)|^pdZ=\int_{\mathbb S^1}|\widehat{U}_k^{(i)}(Z)|^pdZ\ \forall\ p>0$$
 and satisfy simpler recursions as follows:
\

For $k\ge 1,\ i=0,1\ \text{and}\ Z\in\mathbb S^1$, where $(n_{k+1},i)\ne (2,1)$:
  \bul in case $\e_k(0)=0$ we have $s_k(0)=1\ \text{and}\ T_{k+1}-T_k=-(n_{k+1}-1)$ and,
  \begin{align*}
  \widehat{V}_{k+1}^{(i)}(Z)&=Z^{T_{k+1}} \widehat{U}_{k+1}^{(i)}(Z^2)\\ &=
  \sum_{j=1}^{n_{k+1}-1-i}Z^{2(j-1)+T_{k+1}}\widehat{U}_{k}^{(0)}(Z^2)+Z^{2(n_{k+1}-i-1)+T_{k+1}}\widehat{U}_{k}^{(1)}(Z^2)\\ &=
   \sum_{j=1}^{n_{k+1}-1-i}Z^{2(j-1)+T_{k+1}-T_k}\widehat{V}_{k}^{(0)}(Z)+Z^{2(n_{k+1}-i-1)+T_{k+1}-T_k}\widehat{V}_{k}^{(1)}(Z)\\ &=
   \sum_{j=1}^{n_{k+1}-1-i}Z^{2(j-1)-(n_{k+1}-1)}\widehat{V}_{k}^{(0)}(Z)+Z^{n_{k+1}-2i-1}\widehat{V}_{k}^{(1)}(Z)
  \end{align*}

In case $\e_k(0)=1$, we have $T_{k+1}-T_k=T_{k+1}+T_k+2s_k(0)=1$. 
 \
 
\bul in this case, for $n_{k+1}-1-i\in 2\mathbb Z$,
  \begin{align*}
     &\widehat{V}_{k+1}^{(i)}(Z)=Z^{T_{k+1}} \widehat{U}_{k+1}^{(i)}(Z^2)\\ &
=\frac{n_{k+1}-1-i}2 (Z^{T_{k+1}}\widehat{U}_{k}^{(0)}(Z^2) + Z^{2s_k(0)+T_{k+1}}\widehat{U}_{k}^{(0)}(Z^{-2}))+Z^{T_{k+1}}\widehat{U}_{k}^{(1)}(Z^2)\\ &
 =\frac{n_{k+1}-1-i}2 (Z^{T_{k+1}-T_k}\widehat{V}_{k}^{(0)}(Z) + Z^{2s_k(0)+T_{k+1}+T_k}\widehat{V}_{k}^{(0)}(Z^{-1}))+Z^{T_{k+1}-T_k}\widehat{V}_{k}^{(1)}(Z)\\ &
  =\frac{n_{k+1}-1-i}2 (Z\widehat{V}_{k}^{(0)}(Z) + Z\widehat{V}_{k}^{(0)}(Z^{-1}))+Z\widehat{V}_{k}^{(1)}(Z)\end{align*}

  \bul  in this case, for $n_{k+1}-1-i\in 2\mathbb Z+1$,
    \begin{align*}
      &\widehat{V}_{k+1}^{(i)}(Z)=Z^{T_{k+1}} \widehat{U}_{k+1}^{(i)}(Z^2)\\ &=
      \frac{n_{k+1}-i}2 Z^{T_{k+1}} \widehat{U}_{k}^{(0)}(Z^2) + \frac{n_{k+1}-2-i}2Z^{2s_k(0)+T_{k+1}}\widehat{U}_{k}^{(0)}(Z^{-2})+Z^{2s_k(0)+T_{k+1}}\widehat{U}_{k}^{(1)}(Z^{-2})
   \\ &=
      \frac{n_{k+1}-i}2 Z^{T_{k+1}-T_k} \widehat{V}_{k}^{(0)}(Z) + \frac{n_{k+1}-2-i}2Z^{2s_k(0)+T_{k+1}+T_k}\widehat{V}_{k}^{(0)}(Z^{-1})+
      Z^{2s_k(0)+T_{k+1}+T_k}\widehat{V}_{k}^{(1)}(Z^{-1})\\ &=
      \frac{n_{k+1}-i}2 Z \widehat{V}_{k}^{(0)}(Z) + \frac{n_{k+1}-2-i}2Z\widehat{V}_{k}^{(0)}(Z^{-1})+
      Z\widehat{V}_{k}^{(1)}(Z^{-1})
      \end{align*}

 \bul  and
   $\widehat{V}_{k+1}^{(1)}(Z)=Z\widehat{V}_{k}^{(1)}(Z)$ in case $n_{k+1}=2$. 
\section*{\S4   Extracting the {\tt ARW}}   
\subsection*{Associated sequence of temporal probabilities}
\

   Now let $P_k^{(i)}\in\mathcal P(\mathbb Z)$ be defined by 
   $$P_k^{(i)}(\nu):=\frac{V^{(i)}_k(\nu)}{\ell_k(i)}.$$
   We call these ``temporal probabilities'' because
   $$P_k^{(0)}(\nu)=\frac1{\ell_k(0)}\#\{1\le j\le \ell_k(0):\ 2\v_j(0)+T_k=\nu\}.$$
The {\it generating function} of $P_k^{(i)}$ is $\Phi_k^{(i)}:\Bbb S^1\to\Bbb C$ defined by

$$\Phi_k^{(i)}(Z):=\sum_{\nu\in\mathbb Z}P_k^{(i)}(\nu)Z^\nu=\frac{\widehat{V}_k^{(i)}(Z)}{\ell_k(i)}\ \ \ (Z\in\mathbb S^1).$$
   Set 
\begin{align*}
&\phi_N(Z):=\frac1N\sum_{k=0}^{N-1}Z^k\ \ (N\ge 1)\ \text{and}\ \phi_0=\phi_{-1}\equiv 0;\\ &
 p_{k+1}(i):=1-\frac{\ell_k(1)}{\ell_{k+1}(i)}=\frac{(n_{k+1}-1-i)\ell_{k}(0)}{\ell_{k+1}(i)};\\ & 
 q_N:=1-\frac{\lfl \frac{N}2\rfl}N,\ \ \ (N\ge 1)\ \text{and}\ q_0:=0.
\end{align*}

\subsection*{Generating function transitions}
 \
 
 For $k\ge 1,\ i=0,1\ \text{and}\ Z\in\mathbb S^1$  we have, 
 \bul in case $\e_k(0)=0$,
  \begin{align*}
  \Phi_{k+1}^{(i)}(Z)&=
   \frac{\ell_k(0)}{\ell_{k+1}(i)}\sum_{j=1}^{n_{k+1}-1-i}Z^{2(j-1)-(n_{k+1}-1)}\Phi_{k}^{(0)}(Z)+\frac1{\ell_k(i)}Z^{n_{k+1}-2i-1}\Phi_{k}^{(1)}(Z)\\ &=
   p_{k+1}(i)Z^{-(n_{k+1}-1)}\phi_{n_{k+1}-1-i}(Z^2)\Phi_{k}^{(0)}(Z)+(1-p_{k+1}(i))Z^{n_{k+1}-2i-1}\Phi_{k}^{(1)}(Z)
  \end{align*}

 \
 
\bul in case  $\e_k(0)=1\ \text{and}\ 2\mid n_{k+1}-1-i$,
  \begin{align*}
     &\Phi_{k+1}^{(i)}(Z)\\ &=
     \frac{(n_{k+1}-1-i)\ell_k(0)}{2\ell_{k+1}(i)} (Z\Phi_{k}^{(0)}(Z) + Z\Phi_{k}^{(0)}(Z^{-1}))+\frac1{\ell_k(i)}Z\Phi_{k}^{(1)}(Z)\\ &=
     \frac{p_{k+1}(i)}{2} (Z\Phi_{k}^{(0)}(Z) + Z\Phi_{k}^{(0)}(Z^{-1}))+(1-p_{k+1}(i))Z\Phi_{k}^{(1)}(Z)
     \end{align*}

  \bul  in case  $\e_k(0)=1\ \text{and}\ 2\nmid n_{k+1}-1-i$,
    \begin{align*}
      &\Phi_{k+1}^{(i)}(Z)\\ &=
      \frac{(n_{k+1}-i)\ell_k(0)}{2\ell_{k+1}(i)}  Z \Phi_{k}^{(0)}(Z) + \frac{(n_{k+1}-2-i)\ell_k(0)}{2\ell_{k+1}(i)} Z\Phi_{k}^{(0)}(Z^{-1})+
      \frac1{\ell_k(i)}Z\Phi_{k}^{(1)}(Z^{-1})\\ &=
     p_{k+1}(i)q_{n_{k+1}-i-1} Z \Phi_{k}^{(0)}(Z) +  p_{k+1}(i)(1-q_{n_{k+1}-i-1}) Z\Phi_{k}^{(0)}(Z^{-1})+(1-p_{k+1}(i))Z\Phi_{k}^{(1)}(Z^{-1})
      \end{align*}

\subsection*{Smelling the {\tt RAT}s}
\  

Next, for given $\a$ with $2\a=[n_1,n_2,\dots]$ we'll construct  {\tt ARW}s 
$$X^{(k)}=(X^{(k)}(0),X^{(k)}(1))\ \ (k\ge 1)$$ so that
\begin{align*}\tag{\Biohazard}
 P([X^{(k)}(i)=\nu])=P_k^{(i)}(\nu).
\end{align*}
We call an {\tt ARW} satisfying (\Biohazard) an $\a$-{\tt ARW}.

Let
\begin{align*}\tag{\Radioactivity}
&\mathcal N_N\in\text{\tt RV}\,(\mathbb Z),\ \ P(\mathcal N_N=2k-(N-1))=\frac1N,\ 0\le k\le N-1;\\ &
 x_k(i)\in\text{\tt RV}\,(\{0,1\}),\ P(x_k(i)=1)=p_k(i);\\ & 
y_N \in\text{\tt RV}\,(\{0,1\}),\ P(y_N=1)=q_N.\end{align*}
\subsection*{ Flip type {\tt RAT}s}
\

We call a  {\tt RAT} $F=(a,b)\in M_{d\x d}(\{-1,0,1\})\x\mathbb R^d$, {\it of flip type} if
$$a(k,L)=1_{[\frak L(k,a)=L]}\widetilde{a}(k,L)$$
where each $\frak L(k,a)$ is a RV with values in $\{1,2,\dots,d\}$ and each $\widetilde{a}(k,L)$ is a RV with values in $\{-1,1\}$.
\

 Flip type  is preserved under 
composition. For $F'=(a',b')\ \text{and}\ F= (a,b)$,

$$F'\circ F=(a',b')\circ (a,b)=(a'a,a'b+b')$$
and 
$$\frak L(k,a'a)=\frak L(\frak L(k,a'),a).$$

\

The $\a$-{\tt ARW}s to be constructed will be generated by {\tt RAT}s of flip type. Indeed,  henceforward, we only consider flip type {\tt RAT}s

\subsection*{Linear recursion for characteristic functions of flip type {\tt ARW}s}
\

Let $((a^{(n)},b^{(n)}):\ n\ge 1)$ be a flip type {\tt RAT} sequence and consider the generated {\tt ARW} 
$$X^{(0)}=0,\ X^{(n)}=a^{(n)}X^{(n-1)}+b^{(n)}.$$ The characteristic functions of the coordinates of $(X^{(n)},-X^{(n)})$ satisfy a linear recursion.
\

In the special case where $P([a^{(n)}_{k,\ell}=-1])=0\ \forall\ n\ge 1,\ 1\le k,\ell\le d$, there is a simpler linear recursion for the characteristic functions of the coordinates of 
$X^{(n)}$.
\

Writing for the $\mathbb R^d$-valued random variable $X=(X_1,\dots,X_d)$:
$$V_X(\th):=\left(\begin{matrix}  \widehat{\Phi}_{X_1}(\th)\\ \vdots\\ \widehat{\Phi}_{X_d}(\th)
\\  \widehat{\Phi}_{X_1}(-\th)\\ \vdots\\ \widehat{\Phi}_{X_d}(-\th) \end{matrix}\right),$$
where $\widehat{\Phi}_{Y}(\th):=E(e^{i\th Y})$ denotes the characteristic function of the $\mathbb R$-valued random variable $Y$;

and for the flip type {\tt RAT} $(a,b)$ independent of $X$ :
$$ X'=aX+b,$$
we have that
$$V_{X'}(\th)=P(\th)V_X(\th).$$
Here $P(\th)\in M_{2d\x 2d}(\mathbb C)$ is given by
$$P_{k,L}(\th)=P_{k,L}(0)\widehat{\Phi}_{C_{k,L}}(\th)$$ where
$P(0)$ is a stochastic matrix and  $C_{k,L}\ \ (1\le k,L\le 2d)$ are random variables.
\

Specifically:
\begin{align*}& P_{k,L}(0)=P([\frak L(k,a)=L])P([\widetilde{a}_{k,L}=1])\ \text{and}\ \ C_{k,L}=b_{k,L,1};\\ &
P_{k,d+L}(0)=P([\frak L(k,a)=L])P([\widetilde{a}_{k,L}=-1])\ \ \text{and}\ \ C_{k,d+L}=b_{k,L,-1};\\ &
 P_{d+k,L}(0)=P([\frak L(k,a)=L])P([\widetilde{a}_{k,L}=-1])\ \ \text{and}\ \ C_{d+k,L}=-b_{k,L,-1};\\ &
 P_{d+k,d+L}(0)=P([\frak L(k,a)=L])P([\widetilde{a}_{k,L}=1])\ \ \text{and}\ \ C_{d+k,d+L}=-b_{k,L,1}.
\end{align*}
Here, for $1\le k,L\le d,\ J\in\mathbb Z$ and for $\e=\pm 1$, 
\begin{align*}
P([b_{k,L,\e}=J]):=P([b_k=J]|[{a}_{k,L}=\e]).
\end{align*}
Equivalently, for $i,j=0,1\ \text{and}\ \e=1-2i,\ \d=1-2j$
$$ P_{id+k,jd+L}(\th)=P([\frak L(k,a)=L])P([\widetilde{a}_{k,L}=\d])\widehat{\Phi}_{\e\d b_{k,\ell,\d}}(\th).$$
 We'll refer to the function $P:\mathbb R\to M_{2d\x 2d}(\mathbb C)$ as the {\it characteristic function} of the flip type {\tt RAT}:
(abbr.  {\tt RAT-CF}).

 \subsection*{Construction procedures}
 \
 
 Fix a sequence of independent random vectors
 $$(x_{k+1}(i),y_{n_{k+1}-1-i},\mathcal N_{n_{k+1}-1-i}),\ \ i=0,1)_{k\ge 1}$$ whose marginals are determined by $2\a=[n_1,n_2,\dots]$
and (\Radioactivity).  
\

Define $X^{(k)}=(X^{(k)}(0),X^{(k)}(1))\in\text{\tt RV}\,(\mathbb Z^2)$ by $X^{(0)}(i)=0,\ \ i=0,1$ and
 \bul in case $\e_k(0)=0$,
  $$X^{(k+1)}(i)=x_k(i)X^{(k)}(0)+(1-x_k(i))X^{(k)}(1)+x_k(i)\mathcal N_{n_{k+1}-1-i}+(1-x_k(i))(n_{k+1}-2i-1);$$
\bul in case $\e_k(0)=1\ \text{and}\ \ n_{k+1}-1-i\in 2\mathbb Z$,
  $$X^{(k+1)}(i)=x_k(i)(2y_{n_{k+1}-i}-1)X^{(k)}(0)+(1-x_k(i))X^{(k)}(1)+1$$

  \bul in case $\e_k(0)=1\ \text{and}\ \ n_{k+1}-1-i\in 2\mathbb Z+1$,
    $$X^{(k+1)}(i)=x_k(i)(2y_{n_{k+1}-i}-1)X^{(k)}(0)-(1-x_k(i))X^{(k)}(1)+1.$$
    The form of the {\tt RAT}s and the independence of the random vectors implies that $(X^{(k)}:\ k\ge 0)$ is indeed a flip type {\tt ARW}.
\proclaim{Construction Lemma 4.1}\ \ A flip type  {\tt ARW} is an $\a$-{\tt ARW} iff it is defined according to a construction procedure as above .\endproclaim\demo{Proof}
 Evidently,  (\Biohazard) holds if and only if the {\tt ARW} generating function transitions are 
 the same as those established above for the $(\Phi_k^{(i)}(Z):\ k\ge 0,\ i=0,1)$. 
 The latter correspond to linear recursions by {\tt RAT-CF}s of flip type {\tt RAT}s  as above.
 \ \Checkedbox

 \
 
\subsection*{Remark} \ The  plethora of  $\a$-{\tt ARW}s arises  because of  the variety of possible joint distributions of the random vectors
$$\{(x_{k+1}(i),y_{n_{k+1}-1-i},\mathcal N_{n_{k+1}-1-i}):\ \ i=0,1\}$$ for fixed $k\ge 1$. 
\subsection*{$\a$-{\tt RAT} sequences}
\

An {\it $\a$-{\tt RAT} sequence} is a flip type {\tt RAT} sequence  $(F_k=(a^{(k)},b^{(k)}):\ k\ge 1)$  which generates an $\a$-{\tt ARW} $(X^{(k)}:\ k\ge 0)$ by
$$X^{(0)}=0,\ X^{(k)}=a^{(k)}X^{(k-1)}+b^{(k)}.$$ 
\

The  $\a$-{\tt RAT} sequences $(F_k=(a^{(k)},b^{(k)}):\ k\ge 1)$ satisfy, for $n_k\ge 3$, 

\begin{align*}& 
 \tag{0}a^{(k)} =\left(\begin{matrix}  x_k(0) & 1-x_k(0)
\\  x_k(1) & 1-x_k(1) \end{matrix}\right),\ \ \ \ \ \ \ \ \ \ \ \ \e_{k-1}=0;\\ &
\\ &\tag{1}a^{(k)} 
=\left(\begin{matrix}  x_k(0)(2y_{n_{k}}-1)  &-(1-x_k(0))
\\   x_k(1)(2y_{n_{k}-1}-1)  & (1-x_k(1)) \end{matrix}\right)\ \ \ \ \ \ \ \ \ \ \e_{k-1}=1\ \text{and}\ 2\mid n_k;\\ &\tag{1}a^{(k)} 
=\left(\begin{matrix}  x_k(0)(2y_{n_{k}}-1)  &(1-x_k(0))
\\   x_k(1)(2y_{n_{k}-1}-1)  & -(1-x_k(1)) \end{matrix}\right)\ \ \ \ \ \ \ \ \ \e_{k-1}=1\ \text{and}\  2\nmid n_k
\end{align*}
\begin{align*}
\tag{0} b^{(k)} &=\left(\begin{matrix}  x_k(0)\mathcal N_{n_{k}-1}+(1-x_k(0))(n_k-1)
\\  x_k(1)\mathcal N_{n_{k}-2}+(1-x_k(1))(n_k-2) \end{matrix}\right)\ \ \ \ \ \ \ \ \ \ \e_{k-1}=0;
\\ &\tag{1}b^{(k)}
=\left(\begin{matrix}  1
\\ 1  \end{matrix}\right)\ \ \ \ \ \ \ \ \ \ \e_{k-1}=1
\end{align*}
and for $n_{k+1}=2$ by
\begin{align*}& 
 \tag{0}a^{(k)} =\left(\begin{matrix}  x_k(0) & 1-x_k(0)
\\  0 & 1 \end{matrix}\right),\ \ \ \ \ \ \ \ \ \ \ \ \e_{k-1}=0;\\ &
\\ &\tag{1}a^{(k)} 
=\left(\begin{matrix}  x_k(0)(2y_{n_{k}}-1)  &-(1-x_k(0))
\\   0  & 1 \end{matrix}\right)\ \ \ \ \ \ \ \ \ \ \e_{k-1}=1
\end{align*}
\begin{align*}
\tag{0} b^{(k)} &=\left(\begin{matrix}  1-x_k(0)
\\ 0 \end{matrix}\right)\ \ \ \ \ \ \ \ \ \ \e_{k-1}=0;
\\ &\tag{1}b^{(k)}
=\left(\begin{matrix}  1
\\ 1  \end{matrix}\right)\ \ \ \ \ \ \ \ \ \ \e_{k-1}=1.
\end{align*}
\subsection*{Special {\tt RAT}s}\ \ We call {\tt RAT}s of the type defined above  {\it special} ({\tt spec-RAT}s). 
\

We define the {\it parity} of a {\tt spec-RAT} $F$ as above to be {\it even} if it is  defined by equations marked ``(0)'', and to be {\it odd}
if it is  defined by equations marked ``(1)''.  
\

The {\it coefficient} associated to the {\tt spec-RAT}  $F_k$ as above is $n(F_{k}):=n_{k}$.
\

A {\tt spec-RAT} with coefficient $2$ is called {\it trivial}.

  \section*{\S5   Asymptotics of {\tt ARW}s}

\sms{\bf Norm of a matrix}
Throughout this paper, we use the $L^\infty$-operator norm of matrices. Namely, the norm of the matrix $A\in M_{d\x d}(\mathbb C)\cong\mathbb C^{d^2}$ is
$$\|A\|:=\sup\,\{\|Ax\|_\infty:\ x\in\mathbb C^d,\ \|x\|_\infty=1\}=\max_k\,\sum_{\ell=1}^d|A_{k,\ell}|$$
where $\|(x_1,x_2,\dots,x_d)\|_\infty:=\max_{1\le k\le d}\,|x_k|$.
\
\subsection*{Norm of a  {\tt RAT-CF}}
\

Let $P=P_F:\mathbb T\to M_{2d\x 2d}(\mathbb C)$ be the characteristic function of the {\tt RAT} $F=(a,b)$, then for $\th\in\mathbb T$,
\begin{align*}\tag{\Bat}\|&P(\th)\|=\max_{1\le J\le 2d}\sum_{K=1}^{2d}|P_{J,K}(\th)|\\ &=
\max_{1\le k\le d}\,\sum_{L=1}^d\sum_{L=1}^d[P({a}_{k,L}=1)|\widehat{\Phi}_{b_{k,L,1}}(\th)|+P({a}_{k,L}=-1)|\widehat{\Phi}_{b_{k,L,-1}}(\th)|]\\ &\le 1
\end{align*}
with equality iff for 
 some $k,\ b_{k,\ell,\e}$ is a constant random variable $\forall\ 1\le \ell\le d$  and $\e=\pm 1$ with $P({a}_{k,\ell}=\e)>0$.

\

\subsection*{Irreducibility, mean contractivity and balance}
\

We call the {\tt RAT} $F$ {\it irreducible} if $P(a_{k,\ell}(F)\ne 0)>0\ \forall\ 1\le  k,\ell\le d$ and {\it mean contractive} if
$\|E(a(F))\|<1$. Note that $F$ is mean contractive iff for each $1\le k\le d$  $\exists\ 1\le \mu\ne \nu\le d$ so that
$P(a_{k,\mu}\ne 0)>0\ \text{and}\ P(a_{k,\nu}=\e)>0\ \forall\ \e=\pm 1$. 
\

For {\tt spec-RAT}s, mean contractivity entails irreducibility.
\

Nontrivial {\tt spec-RAT}s with odd parity are irreducible and  mean contractive.
\

Nontrivial {\tt spec-RAT}s with even parity are irreducible but not  mean contractive.
\

Trivial  {\tt spec-RAT}s are not irreducible.

\

Call $F$ {\it balanced} if $P(a_{k,\ell}=\e)>0\ \forall\ 1\le k,\ell\le d\ \text{and}\ \e=\pm 1$.
Balance entails both irreducibility and mean contractivity.

If either of $F',\ F$ is irreducible,  mean contractive or balanced, then so is $F'\circ F$.
\

If $F$ is irreducible and $F'$ is mean contractive, then $F'\circ F$ is balanced (but $F\circ F'$ may not be balanced)

\subsection*{Adapted {\tt RAT}s}
\

We'll call the {\tt RAT} $(a,b)\in M_{d\x d}(\mathbb R)\x\mathbb R^d$ 
\bul  {\it adapted} if $\exists\ \th\in\mathbb T$ so that $\|P_{(a,b)}(\th)\|<1$;
\bul    {\it strongly adapted} if  $\|P_{(a,b)}(\th)\|<1\ \forall\ \th\ne 0$;
\bul  {\it partially adapted} if $a$ has  an {\it adapted row} i.e.:  $\exists\ 1\le k\le d\ \text{and}\ \th\in\mathbb T$ so that  
$$\sum_{L=1}^d[P({a}_{k,L}=1)|\widehat{\Phi}_{b_{k,L,1}}(\th)|+P({a}_{k,L}=-1)|\widehat{\Phi}_{b_{k,L,-1}}(\th)|]<1;$$
 equivalently $\exists\ 1\le \ell\le d$ and $\e=\pm 1$ with $P({a}_{k,\ell}=\e)>0,\ b_{k,\ell,\e}$ is a non-constant random variable.  
\

Note that the {\tt RAT} is adapted iff all rows are adapted.

\

The {\tt RAT} $(a,b)\in M_{d\x d}(\mathbb R)\x\mathbb R^d$ with $E(\|b\|^2)<\infty$
 is adapted  iff
\begin{align*}
 \kappa_F & :=\min_k\,-\sum_{\ell=1}^{2d}\frac{d^2}{d\th^2}P_F(\th)_{k,\ell}|_{\th=0}\\ &=
 \min_k\,\sum_{L=1}^dP(\frak L(k,a)=L)[P(\widetilde{a}_{k,L}=1)
 \text{\tt Var}(b_{k,L,1})+P(\widetilde{a}_{k,L}=-1)\text{\tt Var}(b_{k,L,-1})]>0.\end{align*}
 
\

The {\it periodicity group} of the   {\tt RAT} $F=(a,b)$ is
$$\Gamma_F:=\{\th\in\mathbb R:\ \|P_F(\th)\|=1\}.$$
It follows from (\Bat) that
\begin{align*}\tag*{\sun}\G_F=\bigcup_{k=1}^d\bigcap_{1\le\ell\le d}\ \ \bigcap_{\e=\pm 1,P(a_{k,\ell}=\e)>0}\G_{b_{k,\ell,\e}}
\end{align*}

where for a random variable $b$,
\begin{align*}
 \G_b:=\{\th\in\mathbb R:\ |\widehat{\Phi}_b(\th)|=1\}.\end{align*}

This is a closed subgroup of $\mathbb R$ and so, either $\Gamma_F=\mathbb R$ or $F$ is adapted and  $\G_F=g\mathbb Z$ for some $g=g_F>0$.
\

Note that $F$ is strongly adapted iff $\G_F=\{0\}$.
\

In case $F$ is discrete, adapted, $g_F=\frac{2\pi}N$ for some $N\ge 1$ and 
$\G_F=\{\tfrac{2\pi k}N:\ 0\le k\le N-1\}$ is a finite subgroup of $\mathbb T$. Moreover
$$\|P_F(\th+\g)\|=\|P_F(\th)\|\ \forall\ \th\in\mathbb T\ \text{and}\ \g\in\G_F$$
and $\forall\ 0<\e<\tfrac1{2|\G_F|}\ \exists\ \D>0$ so that
$$\|P_F(\th)\|\le 1-\D\ \forall\ \th\in\mathbb T\setminus B(\G_F,\e).$$
Here, for $(X,d)$ a metric space,\ $\G\subset X$ and $\e>0$,
$$B(\G,\e):=\{x\in X:\ \exists\ y\in\G,\ d(x,y)\le\e\}.$$
\sms{\bf Adapted {\tt spec-RAT}s}\ \ Any nontrivial {\tt spec-RAT} $F$  with even parity is:
\sbul adapted if $n(F)\ge 4$ and

\sbul partially adapted if $n(F)=3$.

\

No trivial {\tt spec-RAT}  or {\tt spec-RAT} with odd parity is partially adapted.

\subsection*{Spectral properties of {\tt RAT-CF}s}
\

Let $F$ be a {\tt RAT}.  The {\tt RAT-CF} $P_F:\mathbb R\to M_{2d\x 2d}(\mathbb C)$ is a {\it symmetric function} in the sense that  $P_F(-\th)=\overline{P_F(\th)}$.
\

Noting that
$$P_F(\th)=
\left(\begin{matrix}  Q(\th) & R(\th)  
\\ R(-\th)   & Q(-\th)  \end{matrix}\right)=\left(\begin{matrix}  Q(\th) & R(\th)  
\\ \overline R(\th)   & \overline Q(\th)  \end{matrix}\right)$$
where
$$Q(\th)=\sum_\nu Q_\nu e^{i\nu\th}\ \text{and}\ \ R(\th)=\sum_\nu R_\nu e^{i\nu\th}$$
with each $Q_\nu,\ R_\nu\in M_{d\x d}(\mathbb C)$.

By applying the row and column permutations 
$$(1,2,\dots,2d)\mapsto (d+1,\dots,2d,1,2,\dots,d)$$
(as in \cite{AK}), we see that for $\l\in\mathbb C$,
\begin{align*}
\det(P_F(\th)-\l I)=\det(P_F(-\th)-\l I)=\det(\overline P_F(\th)-\l I).
\end{align*}

It follows  that 
\begin{align*}
 \tag{\Football}\label{football}
\det(P_F(\th)-\l I)=\sum_{k=0}^{2d}c_k(\th)\l^k\  \ \text{where} \ \ c_k:\mathbb T\to\mathbb R\ \text{is even}.
\end{align*}

\subsection*{Spectral theory of {\tt RAT-CF}s}

\

\ \ As seen above a  {\tt RAT-CF} is a perturbation of a positive operator. The spectral theory of perturbations 
(as in chapter III of \cite{HH}) is applicable.
 
 \
 
 Let $F$ be a {\tt RAT} with   $E(\|b(F)\|^2)<\infty$. 
 We call  $P_F(0)$  {\it positive} if $P(0)_{k,\ell}>0\ \forall\ 1\le k,\ell\le d$ 
 (equivalently, $F$ is balanced), and {\it aperiodic} if
 $\exists\ N\ge 1$ with $P_F(0)^N$ is positive.
 \
 
 Suppose that  $P_F(0)$ is aperiodic.
 
 \
 
 By the Perron-Frobenius theorem,
$1=\|P_F(0)\|$ is a simple, dominant eigenvalue (i.e. its eigenspace $\mathbb C\cdot\mathbb{1}$ is one-dimensional and all other eigenvalues are smaller in absolute value). 
\

Since $\th\mapsto P_F(\th)$ is $C^2$, by the implicit function theorem, $\exists\ \e>0$ so that for $|\th|<\e$, 
$P_F(\th)$ also has a  simple, dominant eigenvalue $\l_F(\th)\in \mathbb C$ where   $\l_F:(-\e,\e)\to\mathbb R$ is differentiable.

By (\Football) 
 $\l_F:(-\e,\e)\to\mathbb R$ is a real valued, even function  with  $\l_F'(0)=0$ and $\g_F:=-\l_F^{\pprime}(0)\ge 0$ with equality iff 
 $\l_F\equiv 1$. 
 \
 
 If $F$ is adapted, then $\g_F>0$. However, it may be that $\g_F>0$ and $\|P_F\|\equiv 1$.
 
 \
 
Fix $\pi=\pi_F\in\mathbb R_+^d$ satisfying
$$ \<\pi,\mathbb{1}\>=1\ \text{and}\ P_F(0)^t\pi=\pi.$$
There is a unique eigenvalue function $\eta=\eta_F:(-\e,\e)\to\mathbb C^{2d}$ so that
$$P_F(\th)\eta(\th)=\l(\th)\eta(\th)\ \text{and}\ \<\pi,\eta(\th)\>=1\ \ \forall\ |\th|<\e.$$
It follows that $\eta(0)=\mathbb{1}\ \text{and}\ \eta(-\th)=\overline{\eta(\th)}$.

\proclaim{Theorem 5.1\ Coordinate distributional limits for stationary {\tt ARW}s}
\

Let $F_n=(a_n,b_n)\in M_{d\x d}(\mathbb R)\x\mathbb R^d\ \ (n\ge 1)$ be an iid {\tt RAT} sequence
with each $F_n\overset{\text{\tiny dist}}= F$ with $F$ adapted and $P_{F}(0)$  aperiodic and let
$(X_1,X_2,\dots)$ be the associated {\tt ARW}, then for each $1\le k\le d$ 
\begin{align*}\tag{\tt CLT} \frac{(X_n)_k}{\sqrt n}\xrightarrow[n\to\infty]{\text{\tt\tiny dist}}\mathcal N(0,\g_F)
\end{align*}
If, in addition, $F$ is strongly adapted, then for any bounded interval $J\subset\Bbb R$, 
\begin{align*}\tag{\tt LLT} 
 \sqrt nP([(X_n)_k\in t_n+J])\xrightarrow[n\to\infty,\ \frac{t_n}{\sqrt n}\to x]{} \frac1{2\pi\g_F}e^{-\frac{x^2}{2\g_F}}.
\end{align*}
\endproclaim\demo{Proof sketch}\ \ Let $P=P_F$ and let $\pi\in\mathbb R_+^{2d}$ be so that
$P_F(0)^t\pi=\pi$. For each $\th\in (-\e,\e),\ \exists\ !\ \eta(\th)\in\mathbb C^d$ so that
\begin{align*}P_F(\th)\eta(\th)=\l_i(\th)\eta(\th)\ \text{and}\ \<\eta(\th),\pi\>=1.
\end{align*}
It follows that $\eta:(-\e,\e)\to\mathbb C^d$ is a smooth, symmetric function.
\

Again, for each $\th\in (-\e,\e),\ \exists\ !\ \pi(\th)\in\mathbb C^d$ so that
 for $\th\in (-\e,\e)$,
\begin{align*} P_F(\th)\pi(\th)=\l(\th)\pi(\th)\ \text{and}\  \<\pi(\th),\eta(\th)\>=1.
\end{align*}
It follows that $\pi:(-\e,\e)\to\mathbb C^d$ is also a smooth symmetric function.

Define the projection $N(\th):\mathbb C^d\to\mathbb C^d$ by
$$N(\th)x:=\<x,\pi(\th)\>\eta(\th).$$
We see that 
$$P_F(\th)N(\th)=\l(\th)N(\th),$$
and that $Q(\th):=P_F(\th)-P_F(\th)N(\th):\mathbb C^d\to\mathbb C^d$ satisfies $QN=NQ=0$.
The spectral radius of $Q(0)$,\ $r(Q(0))<1$. By  continuity  of $\th\mapsto Q(\th)$,  for possibly smaller $\e>0$, 
$$r(Q(\th))\le\rho<1\ \forall\ |\th|<\e.$$
Thus
$$V_{X_n}(\th)=P_F(\th)^n\mathbb{1}=\l(\th)^n\eta(\th)+O(\rho^n).$$
The first statement follows directly from this. 
\

The second follows  also since if $F$ is strongly adapted,
$\sup_{|\th|\in[\e,M]}|\l(\th)|<1\ \forall\ 0<\e<M<\infty$. This is seen via standard proofs of the local limit theorem (see \cite{Breiman}).\ \Checkedbox
\sms{\bf Remark:\ Irreducible, positive {\tt RAT}s}
\

\sms  The above does not include the important case of an irreducible {\tt RAT} $(a,b)$ with $P(a_{l,\ell}=-1)=0\ \forall\ k,\ \ell$. 
Here there is possible linear drift, but the d-dimensional vectors
$(\widehat{\Phi}_{X_1^{(n)}},\widehat{\Phi}_{X_2^{(n)}},\dots,\widehat{\Phi}_{X_d^{(n)}})$ satisfy a linear recursion and an analogous spectral 
argument applies.  

\

\sms{\bf Remark:\  Quadratic $\a$.}
\

The following is a special case of a subsequence version of theorem 1.1 in \cite{Beck1} (see also \cite{ADDS}).

\proclaim{Proposition 5.2}\ \ Suppose that $\a\in (0,1)$ is quadratic. 
There are $K,\ L\in\Bbb N$ and constants $c>0,\ \mu\in\mathbb R$ so that for intervals $I\subset\mathbb R$
$$\frac1{\ell_{K+Ln}(0)}\#\left\{1\le j\le \ell_{K+Ln}(0):\ \frac{\v_j^{(\a)}(0)-\mu n}{c\sqrt n}\in\ I\right\}\xrightarrow[n\to\infty]{}\frac1{2\pi}\int_Ie^{\frac{-t^2}2}dt.$$
\endproclaim
Here, $\v_j^{(\a)}(0):=\sum_{t-0}^{j-1}\v\circ r_\a(0)$.
This does not follow directly from {\tt CLT} in theorem 5.1 as the $\a$-{\tt RAT} sequence  is not eventually periodic 
(although it converges to an  eventually periodic {\tt RAT} sequence). 
However,  the transitions of the generating functions of the simplified visit distributions are 
eventually periodic and the result can be deduced from this using arguments in section 5 of  \cite{AK} as was (\dsrailways) 
(see page \pageref{choochoo}) in this case.
\

For general  $\a\in\text{\tt BAD}$ there need be no eventual periodicity of  renormalization related sequences.  New tools are needed. 
 
\section*{\S6  The weak, rough local limit theorem}

\subsection*{Variance of {\tt ARW}s}

\
 
 We call the {\tt RAT} $(a,b)$ {\it centered} if  $E(b)=0$ $\text{and}$ the  {\tt RAT} sequence $({F}_n=(a^{(n)},{b}^{(n)}):\ n\ge 1)$ {\it centered} if each
 {\tt RAT} ${F}_n=(a^{(n)},{b}^{(n)}) $ is centered. 
 \
 
 From the formula
 $$X^{(n)}=F_n(X^{(n-1)})=a^{(n)}X^{(n-1)}+b^{(n)}$$
 we see that
\begin{align*}E(X^{(n)})=E(a^{(n)})E(X^{(n-1)})+E(b^{(n)})\end{align*}
  whence
the {\tt RAT} sequence $({F}_n=(a^{(n)},{b}^{(n)}):\ n\ge 1)$ is  centered if and only if
$$E(X^{(n)})=0\ \forall\ n\ge 1.$$

            \proclaim{6.1 Variance Theorem}
\

Let $(F_n:=(a^{(n)},b^{(n)}): \  n\ge 1)$ be a centered, independent  {\tt RAT} sequence and let
$$X^{(n)}:=F_1^n(0)=F_n\circ F_{n-1}\dots\circ F_1(0)$$
be the corresponding {\tt ARW}, then for $1\le k\le d\ \text{and}\ n\ge 1$,
$$\sum_{\nu=1}^n\min_{1\le L\le d}E([b^{(\nu)}_L]^2)\le E([X^{(n)}_k]^2)\le\sum_{\nu=1}^n\max_{1\le L\le d}E([b^{(\nu)}_L]^2).$$ \endproclaim

\demo{Proof}
\

Set 

\begin{align*}Y^{(n,\nu)} :=a_{\nu+1}^nb^{(\nu)}\ \ (1\le \nu\le n)
\end{align*}
where
 \begin{align*}a_k^N:=\begin{cases}& a^{(N)}a^{(N-1)}\cdots a^{(k)}\ \ \ \ k\le N\\ & \text{Id}\ \ \ \ k>N.\end{cases}  
 \end{align*}
so that
\begin{align*}X^{(n)}=\sum_{\nu=1}^nY^{(n,\nu)}.
\end{align*}

 As above, $E(Y_K^{(n,\nu)})=0$.
 \

 \Par\ \ \  for each $1\le K\le d,\ n\ge 1$, the random variables $\{Y_K^{(n,\nu)}:\ 1\le\nu\le n\}$ are  orthogonal, i.e. 
 $$E(Y_K^{(n,\nu)}Y_K^{(n,\mu)})=0 \ \forall\ \mu\ne\nu.$$
 \demo{Proof of \P}
\

We have for $1\le K\le d$, that
\begin{align*}\tag{\dsagricultural}
 Y^{(n,\nu)}_K=(a_k^N)_{K,\ell(K,a_{\nu+1}^n)}b^{(\nu)}_{\ell(K,a_{\nu+1}^n)}
\end{align*}

and for $1\le \nu<\mu\le n$, that
\begin{align*}
 &\ell(K,a_{\nu+1}^n)=\ell(\ell(K,a_{\mu+1}^n),a_{\nu+1}^\mu)\ \text{and}\ \\ &
 (a_{\nu+1}^n)_{K,\ell(K,a_{\nu+1}^n)}=(a_{\mu+1}^n)_{K,\ell(K,a_{\mu+1}^n)}(a_{\nu+1}^\mu)_{\ell(K,a_{\mu+1}^n),\ell(\ell(K,a_{\mu+1}^n),a_{\nu+1}^\mu)} 
\end{align*}
 Consequently
\begin{align*}
 Y^{(n,\nu)}_K&Y^{(n,\mu)}_K\\ &=(a_{\nu+1}^n)_{K,\ell(K,a_{\nu+1}^n)}b^{(\nu)}_{\ell(K,a_{\nu+1}^n)}(a_{\mu+1}^n)_{K,\ell(K,a_{\mu+1}^n)}b^{(\mu)}_{\ell(K,a_{\mu+1}^n)}
 \\ &=
 (a_{\nu+1}^\mu)_{\ell(K,a_{\mu+1}^n),\ell(\ell(K,a_{\mu+1}^n),a_{\nu+1}^\mu)}
 b^{(\nu)}(\ell(\ell(K,a_{\mu+1}^n),a_{\nu+1}^\mu))b^{(\mu)}_{\ell(K,a_{\mu+1}^n)}
 \\ &=
 \sum_{L=1}^d1_{[\ell(K,a_{\mu+1}^n)=L]}b^{(\mu)}_L\sum_{M=1}^d1_{[\ell(L,a_{\nu+1}^\mu)=M]}
 (a_{\nu+1}^\mu)_{L,M} b^{(\nu)}_M
\end{align*}
whence by independence and centering,

 \begin{align*}
 E(Y^{(n,\nu)}&_KY^{(n,\mu)}_K)\\ &=
 \sum_{L=1}^dP([\ell(K,a_{\mu+1}^n)=L])E(b^{(\mu)}_L)\sum_{M=1}^dE(1_{[\ell(L,a_{\nu+1}^\mu)=M]}(a_{\nu+1}^\mu)_{L,M})E(
 b^{(\nu)}_M)\\ &=0.\ \ \CheckedBox\P
 \end{align*}
It  follows from (\dsagricultural) that

 \begin{align*}
 E(Y^{(n,\nu)2}_K)=\sum_{L=1}^dP([\ell(K,a_{\nu+1}^n)=L])E(b^{(\nu)2}_L)
 \end{align*}
 and from the above that
 $$E(X^{(n)2}_K)=\sum_{\nu=1}^nE(Y^{(n,\nu)2}_K)).$$
 The Variance Theorem follows from this.\ \ \Checkedbox

\subsection*{Compactness  properties of {\tt RAT} sequences}
\

\

\

We'll say that the {\tt RAT} sequence 
$$(F_n=(a^{(n)},b^{(n)}):\ n\ge 1)\in\text{\tt RV}(M_{d\x d}(\mathbb R)\x\mathbb R^d)^\mathbb N$$

\bul is {\it adapted} if $\exists$ a finite subgroup $\G\le\mathbb T$ (called the {\it adaptivity group})  so that 
\sms $|\widehat{\Phi}_{b_{k,\ell,\pm 1}^{(n)}}(\th)|=1\ \Lra\ \th\in\G$, and
\sms $\forall\ \e>0\ \exists\ \d>0$ so that 
$$|\widehat{\Phi}_{b_{k,\ell,\pm 1}^{(n)}}(\th)|\le 1-\d\ \forall\ \ 1\le k,\ell\le d,\ n\ge 1\ \ \text{and} \ \th\in\mathbb T\setminus B(\G,\d).$$

\bul is {\it uniform} if $\sup\{|E(b_{k,\ell,\e}^{(n)})|:\ <\infty\  \forall\ \ 1\le k,\ell\le d,\ n\ge 1,\ \e=\pm 1$ and
\sms (a) the collection

$$\mathcal F_F:=\{|\overline{b}_{k,\ell,\e}^{(n)}|^2:\ \ 1\le k,\ell\le d,\ n\ge 1,\ \e=\pm 1\}$$
is uniformly integrable   where $\overline{Y}:=Y-E(Y)$ denotes the centering of the integrable random variable $Y$; and 
\sms (b) $\exists\ M>1$ so that
\begin{align*}
 \text{\tt Var}\,(b^{(n)}_{k,\ell,\e})=M^{\pm 1}\sum_{j=1}^d&\text{\tt Var}\,(b^{(n)}_j)\\ &\forall\ n\ge 1,\ 1\le k,\ell\le d,\ \e=\pm 1\ \text{with} \ \ P(a^{(n)}_{k,\ell}=\e)>0.
\end{align*}

\

It is standard to show that uniformity of $((a^{(n)},b^{(n)}):\ n\ge 1)$ implies that $\exists\ M>1$ so that
$\text{\tt Var}\,(b_{k,\ell,\e}^{(n)})=M^{\pm 1}\ \ \forall\ n\ge 1,\ 1\le k,\ell\le d,\ \e=\pm 1$  with  $P(a^{(n)}_{k,\ell}=\e)>0$.

\

The above properties do not entail discreteness and we'll need the {\tt  WRLLT} for certain non-discrete {\tt ARW}s.
\

\proclaim{Theorem  6.2\ \ ({\tt WRLLT})}\ \ Let 
 $(F_n=(a^{(n)},b^{(n)}):\ n\ge 1)\in (M_{d\x d}(\mathbb R)\x\mathbb R^d)^\mathbb N$ be a {\tt RAT} sequence.
 \
 
 Suppose that $(F_n:\ n\ge 1)$ is adapted, centered and uniform and let 
 $(X^{(n)})_{n\ge 1}$ be the corresponding {\tt ARW}.
 \
 
For each $1\le k\le d$ and $1\le p\le 2$,
$$\int_\mathbb T|\widehat{\Phi}_{X^{(n)}_k}(\th)|^pd\th\ \text{\Large $\asymp$}\ \frac1{\sqrt{n}}.$$\endproclaim

\demo{Proof of the  {\tt WRLLT}}\ \ Since the integral decreases with $p$, it suffices to show that $\exists\ M>1$ so that $\forall\ n\in\mathbb N$ large,
$$\text{\rm (a)}\ \ \ \ \int_\mathbb T|\widehat{\Phi}_{X^{(n)}_k}(\th)|^2d\th\ \ge\ \frac1{M\sqrt{n}}\ \ \text{and}\ \ \ 
\text{\rm (b)}\ \ \ \ \int_\mathbb T|\widehat{\Phi}_{X^{(n)}_k}(\th)|d\th\ \le\ \frac{M}{\sqrt{n}}.$$
\

\demo{Proof of (a)}\ \ It follows from adaptedness, centeredness and the Variance theorem that $\exists\ \G>0$ so that
$$E(X^{(n)2}_k)\le \G n\ \ \ (1\le k\le d).$$
Next, fix $M=2\sqrt{\G}$, then by Chebyshev's inequality,

$$P([|X^{(n)}_k|\le M\sqrt n])\ge\frac34.$$
Now fix $\D>0$ so that 
$$|1-e^{it}|<\frac14\ \forall\ |t|<\D.$$
We have
\begin{align*}\sqrt n\int_\mathbb T|E(e^{i\th X^{(n)}_k})|^2d\th &=\int_{-\pi\sqrt n}^{\pi\sqrt n}|E(\exp[it\frac{X^{(n)}_k}{\sqrt n})|^2dt\\ &\ge
\int_{-\frac{\D}M}^{\frac{\D}M}|E(\exp[it\frac{X^{(n)}_k}{\sqrt n}])|^2dt.
 \end{align*}
For $|t|<\frac{\D}M$, we have
\begin{align*}
 |E(\exp[it\frac{X^{(n)}_k}{\sqrt n})|&\ge |E(\exp[it\frac{X^{(n)}_k}{\sqrt n}]1_{[|X^{(n)}_k|<M\sqrt n]})|-P([|X^{(n)}_k|\ge M\sqrt n]\\ &\ge 
 \frac34\cdot\frac34-\frac14\\ &=\frac{5}{16}\end{align*}
 whence
 \begin{align*}\sqrt n\int_\mathbb T|E(e^{i\th X^{(n)}_k})|^2d\th &\ge
\int_{-\frac{\D}M}^{\frac{\D}M}|E(\exp[it\frac{X^{(n)}_k}{\sqrt n}])|^2dt\\ &\ge \frac{2\D}M\cdot\frac{25}{256}.\ \ \ \CheckedBox\text{\rm (a)}
 \end{align*}
\demo{Proof of (b)}\ \ Suppose that $\G=\{\tfrac{k}N:\ 0\le k\le N-1\}\subset\mathbb T$. 
\

We claim first that 
\Par \ $\exists\ \rho\in (0,1),\ 0<\D<\frac1{2N}$ and $\b>0$  so that for each $n\ge 1$ and $\g\in \G_{F_n}$,
\begin{align*}
 &\tag{i}\|P_{F_n}(\g+\th)\|\le 1-\b\th^2\ \ \forall\ |\th|<\D;\ \text{and}\\ &
 \tag{ii}\|P_{F_n}(\g+\th)\|\le \rho\ \ \forall\ \th\notin B(\F,\D)\\ &
 \tag{iii}\rho\le 1-\b\D^2.\end{align*}
 
 \demo{Proof of {\rm \P(i)}}\ \ Let $Y$ be a random variable with finite second moment and let $W=Y-Y'$ be its symmetrization.
 We have for $\th\in\mathbb T\cong (-\frac12,\frac12)$ that
 $$|\Phi_Y(\th)|^2=\Phi_W(\th)=1-2E(\sin^2(\tfrac{W\th}2)).$$
 Now
 \begin{align*}E(\sin^2(\tfrac{W\th}2))&\ge E(1_{[|W|\le\frac{\pi}\th]}\sin^2(\tfrac{W\th}2))\\ &\ge
 \frac{2\th^2}{\pi^2}E(1_{[|W|\le\frac{\pi}\th]}W^2)\\ &=
 \frac{2\th^2}{\pi^2}E(W^2)(1-\D_Y(\th))
 \end{align*}
where
$$\D_Y(\th):=\frac{E(1_{[|W|>\frac{\pi}\th]}W^2)}{E(W^2)}.$$
It follows from uniformity that $\exists\ \D>0,\ \D<\tfrac1{2N}$ so that
$$\D_{\overline{b}_{k,\ell,\e}^{(n)}}<\frac12\ \ \forall\ |\th|<\D,\ 1\le k,\ell\le d,\ n\ge 1$$
whence for $|\th|<\D\ \text{and}\ \g\in\G_{F_n}$,
 \begin{align*}
  \|P_{F_n}(\g&+\th)\|=\|P_{F_n}(\th)\|\\ &=
  \max_{1\le k\le d}\,\sum_{L=1}^dP(\frak L(k,a)=L)[P(\widetilde{a}_{k,L}=1)|\widehat{\Phi}_{b_{k,L,1}}(\th)|+P(\widetilde{a}_{k,L}=-1)|\widehat{\Phi}_{b_{k,L,-1}}(\th)|]\\ &  
  \le  1-\frac{\th^2}{\pi^2}\min_{1\le k\le d}\,\,\sum_{L=1}^dP(\frak L(k,a)=L)[P(\widetilde{a}_{k,L}=1)\text{\tt Var}(b^{(n)}_{k,\ell,1})
  +P(\widetilde{a}_{k,L}=-1)\text{\tt Var}(b^{(n)}_{k,\ell,-1})]\\ &
   = 1-\frac{\kappa_{F_n}\th^2}{\pi^2}\\ &\le 1-\b\th^2
 \end{align*}
 for any $0<\b\le \min_{n\ge 1}\frac{\kappa_{F_n}}{\pi^2}$ which latter is positive by     uniformity. \ \Checkedbox(i)
\demo{Proof of {\rm (ii)\   $\text{and}$ \ (iii)}}
Statement (ii) follows from  adaptedness. Statement (iii) can be obtained by shrinking $\b$.\ \Checkedbox (ii), (iii)\ $\text{and}$\ \P.

To complete the proof of (b), we have
\begin{align*}
 |\widehat{\Phi}_{X^{(n)}_k}(\th)|&=|(P_{F_n}(\th)P_{F_{n-1}}(\th)\cdots P_{F_1}(\th)\mathbb{1})_k|\\ &
 \le \prod_{k=1}^n\|P_k(\th)\|.
\end{align*}
By \P(i), for $n\ge 1,\ \g\in\G_{F_n}$ we have
$$\|P_{F_n}(\g+\th)\|\le 1-\b \th^2\ \ \forall\ |\th|<\D$$ and by \P(ii)
$$\|P_{F_n}(t)\|\le \rho\ \ \forall\ t\notin B(\G_{F_n},\D).$$
It follows that for $\g\in\G\setminus B(\G_{F_n},\D)$ and $|\th|<\D$, we have
\begin{align*}\|P_{F_n}(\g+\th)\|&\le\rho \ \ \text{by \P(ii)}\\ &\le
1-\b\th^2 \ \ \text{by \P(iii)}.
 \end{align*}

Thus for $\g\in\G$:
$$ |\widehat{\Phi}_{X^{(n)}_k}(\g+\th)|\le (1-\b\th^2)^n\ \ \ \ |\th|<\D$$
and
$$ |\widehat{\Phi}_{X^{(n)}_k}(t)|\le
\rho^n \ \ \forall\ t\notin B(\G,\D),$$
whence
\begin{align*}\int_\mathbb T|\widehat{\Phi}_{X^{(n)}_k}(\th)|d\th &\le 
(\int_{B(\G,\D)}+\int_{\mathbb T\setminus B(\G,\D)}) \prod_{k=1}^n\|P_k(\th)\|d\th\\ &\le
\sum_{\g\in\G}\int_{-\D}^{\D}\prod_{k=1}^n\|P_k(\g+\th)\|d\th+\int_{\mathbb T\setminus B(\G,\D)}\prod_{k=1}^n\|P_k(\th)\|d\th\\ &\le
N\int_{-\D}^{\D}(1-\b\th^2)^nd\th+\rho^n\\ &
\le\frac{N}{\sqrt n}\int_{-\D\sqrt{n}}^{\D\sqrt{n}}(1-\frac{\b t^2}n)^ndt+\rho^n\\ &\le
\frac{N}{\sqrt n}\int_\mathbb Re^{-\b t^2}dt  +\rho^n\\ &\propto \frac1{\sqrt n}.\ \ \ \CheckedBox\text{(b)}\ \text{and}\ \text{\tt WRLLT} \end{align*}

\

\

\

\section*{\S7 Centering of a {\tt RAT} sequence}

\

Let 
 $(F_n=(a^{(n)},b^{(n)}):\ n\ge 1)$ be an independent {\tt RAT} sequence  and let $X^{(n)}:=F_n\circ F_{n-1}\circ\dots\circ F_1(0)$ $\text{and}\ c_n:=E(X^{(n)})$. Set  $\widetilde{X}^{(n)}=X^{(n)}-c_n$, then
 \begin{align*}\widetilde{X}^{(n)}&=X^{(n)}-c_n\\ &=a^{(n)}X^{(n-1)}+b^{(n)}-c_n
 \\ &=a^{(n)}\widetilde{X}^{(n-1)}+b^{(n)}-c_n+a^{(n)}c_{n-1}\\ &=:\widetilde{F_n}(\widetilde{X}^{(n-1)})
 \\ &=\widetilde{F}_n\circ \widetilde{F}_{n-1}\circ\dots\circ\widetilde{F}_1(0).\end{align*} 
 and so $(\widetilde{X}^{(n)})_{n\ge 1}$ is a centered {\tt ARW} with the corresponding centered independent,  {\tt RAT} sequence
 $(\widetilde{F}_n=(a^{(n)},\widetilde{b}^{(n)}):\ n\ge 1)$ where
  $$\widetilde{b}^{(n)}=b^{(n)}-c_n+a^{(n)}c_{n-1}=0.$$ 

The {\tt ARW} $\text{and}$ associated {\tt RAT} sequence above are unique. We call them the {\it centerings} of the {\tt ARW}  
$\text{and}$ associated {\tt RAT} sequence respectively.
\proclaim{Proposition 7.1}\ \ \ Let $(X^{(n)}:\ n\ge 1)$ be an {\tt ARW} with centering $(\widetilde{X}^{(n)})_{n\ge 1}$, then for each $n\ge 1\ \text{and}\ 1\le k\le d$,
$$\text{\tt Var}(X_k^{(n)})=E(\widetilde{X}_k^{(n)2}).$$\endproclaim
\subsection*{Adaptedness preservation}
\

The centering of a discrete {\tt ARW} may not be discrete, but centering does not affect adaptedness or any of the other compactness properties. 
\

This is because if $(F_n:\ n\ge 1)$ is a {\tt RAT} with centering $(\widetilde{F}_n:\ n\ge 1)$, then
$$\widetilde{F}_n=(a(F_n),b(F)-c^{(n)}+a(F)c^{(n-1)})$$ where $c^{(n)}\in\mathbb R^d\ n\ge 1$ are constant. 
\

Thus, for
$1\le k,\ell\le d\ \text{and}\ \e=\pm 1$ with $P(a_{k,\ell}=\e)>0$, the random variable
 $\widetilde{b}_{k,\ell,\e}$ is a constant translation of ${b}_{k,\ell,\e} $:
  $$\widetilde{b}_{k,\ell,\e}(\widetilde{F}_n)={b}_{k,\ell,\e}(F_n)-c_k^{(n)}+\e c_\ell^{(n-1)}.$$
  It follows from \sun that the adaptedness of $(F_n:\ n\ge 1)$ is equivalent to that of $(\widetilde{F}_n:\ n\ge 1)$ and the adaptivity groups are the same.
  
  \

We'll need conditions for $\sup_{n\ge 1}\|b^{(n)}-\widetilde{b}^{(n)}\|<\infty$. a.s..
\subsection*{Bounded mean fluctuation}
\

We say that the {\tt ARW} $(X^{(n)}:=F_n\circ F_{n-1}\circ\dots\circ F_1(0):\ n\ge 1)$ (and its associated {\tt RAT} sequence) has {\it bounded mean fluctuations} {\tt BMF} if
$$\sup_{n\ge 1}\|E(X^{(n)})\|<\infty.$$
\proclaim{7.2  Bounded Centering Proposition} \ \ Let $(F_n=(a^{(n)},{b^{(n)}}):\ n\ge 1)$ be a {\tt RAT} sequence with centering
$(\widetilde{F_n}=(a^{(n)},\widetilde{b^{(n)}}):\ n\ge 1)$.
\

If $(F_n=(a^{(n)},b^{(n)}):\ n\ge 1)$ has {\tt BMF}, then $\exists\ M>1$ so that
$$\|\widetilde{b^{(n)}}-b^{(n)}\|\le M\ \ \text{a.s.}\ \ \forall\ n\ge 1.$$\endproclaim
\demo{Proof}\ \ Since
 $$\widetilde{b}^{(n)}-b^{(n)}-c_n+a^{(n)}c_{n-1}$$ and $\|a^{(n)}\|\le 1$ a.s., the proposition holds with 
$$M=2\sup_{n\ge 1}\|E(X^{(n)}\|.\ \ \ \CheckedBox$$
The following gives sufficient conditions for {\tt BMF} of a {\tt RAT} sequence in terms of uniform mean boundedness and uniform mean contractivity.
\proclaim{7.3 Lemma } \ \ Let $(F_n=(a^{(n)},{b^{(n)}}):\ n\ge 1)$ be a {\tt RAT} sequence.  If  $\exists\ M>0\ \text{and}\ \rho\in (0,1)$ so that
$$\|E(b^{(k)})\|\le M\ \text{and}\ \|E(a^{(k)})\|\le\rho\ \forall\ k\ge 1,$$ then $(F_n=(a^{(n)},{b^{(n)}}):\ n\ge 1)$ has {\tt BMF}.\endproclaim\demo{Proof}
\begin{align*}X^{(n)}:=F_n\circ F_{n-1}\dots\circ F_1(0)=\sum_{\nu=1}^na_{\nu+1}^nb^{(\nu)}
\end{align*} 
where
 \begin{align*}a_k^N:=\begin{cases}& a^{(N)}a^{(N-1)}\cdots a^{(k)}\ \ \ \ k\le N\\ & \text{Id}\ \ \ \ k>N,\end{cases}  
 \end{align*}

we have by independence that
$$E(X^{(n)})=\sum_{\nu=1}^nE(a_{\nu+1}^n)E(b^{(\nu)})=\sum_{\nu=1}^n\prod_{k=\nu+1}^nE(a^{(k)})E(b^{(\nu)})$$
whence
$$\|E(X^{(n)})\|\le\sum_{\nu=1}^n\prod_{k=\nu+1}^n\|E(a^{(k)})\|\|E(b^{(\nu)})\|\le M\sum_{k=1}^\infty\rho^k=\frac{M\rho}{1-\rho}<\infty.
\ \ \ \CheckedBox$$
\

 \section*{\S8  Proof of the main result}
 \
 
 Our first step in the proof of (\dsrailways) (as on page \pageref{choochoo}) is to show that for $\a\in\text{\tt BAD}$, each  $\a$-{\tt ARW} satisfies the 
 {\tt WRLLT} along a syndetic subsequence.
 \subsection*{A canonical subsequence} \  Define $\nu:\mathbb N_2^\mathbb N\to\mathbb N\cup\{\infty\}$ by
$$\nu(n_1,n_2,\dots):=\min\,\{J\ge 1:\ \sum_{k=1}^J1_{[n_k>2]}\ge 4\}.$$

The number $\b=2\a=[n_1,n_2,\dots]$ is irrational iff $\#\{k\ge 1:\ n_k>2\}=\infty$ so $\nu(n_{k+1},n_{k+2},\dots)<\infty\ \ \forall\ k\ge 1.$ 
\

Define
$$\nu_0=0,\ \nu_{k+1}:=\nu_k+\nu(n_{{\nu_k}+1},n_{{\nu_k}+2}\dots),$$
then  $\nu_k<\infty\ \forall\ k\ge 1$.

By \cite{KN},  $\a\in\text{\tt BAD}$   iff
\sms (i) $\sup_{k\ge 1}\,n_k<\infty$ and (ii)\ $\sup\,\{J\ge 1:\ \exists\ k,\ n_{k+j}=2\ \forall\ 1\le j\le J\}<\infty$.

Thus  in case  $\a\in\text{\tt BAD}$,\ \ \ $\sup_{k\ge 1}(\nu_{k+1}-\nu_k)<\infty$ i.e. $\{\nu_k\}_{k\ge 1}$ is syndetic..
\

\subsection*{The grouping  theorem}\ \ Let $(F_n=(a^{(n)},b^{(n)}):\ n\ge 1)$ be a  $\a$-{\tt RAT} sequence. The {\it canonical grouping} of $\a$ 
(or $(F_n:\ n\ge 1)$)
is the {\tt RAT} sequence  $(G_n:\ n\ge 1)$ defined by
$$G_0=\text{\tt Id},\ G_{k+1}:=F_{\nu_{k+1}}\circ F_{\nu_{k+1}-1}\circ\dots\circ F_{\nu_k+1}.$$ 
\

\proclaim{8.1  Grouping theorem} \ \ Let $\a\in\text{\tt BAD}$,  let $(F_k:\ k\ge 1)$ be the associated  {\tt RAT}, then the canonical grouping  $(G_k:\ k\ge 1)$ 
 of  $(F_k:\ k\ge 1)$ has {\tt BMF} and  is 
 adapted.\endproclaim
 
 The proof of the Grouping theorem proceeds via:

\subsection*{Compact {\tt RAT} collections}\label{compact}
\

\ \ We'll call the collection $\frak F\subset \text{\tt RV}\,(M_{d\x d}(\mathbb Z^d)\x\mathbb Z^d)$ of flip type {\tt RAT}s
 {\it compact} if 
\sms $\exists\ M=M_\frak F>1$ so that $\forall\ F:=(a,b)\in\mathfrak F$ and

  $1\le k,\ell\le d,\ \e=0,\pm 1,\ \nu\in\mathbb Z$ we have
\begin{align*} &\tag{i} P([a_{k,\ell}=\e]\cap [b_k=\nu])>0\ \Lra\ P([a_{k,\ell}=\e]\cap [b_k=\nu])\ge\tfrac1M\ \text{and}\\ &\tag{ii}
 \|b\|_\infty:=\sup_{P(b=\b)>0}\,\|\b\|\ \le \ M.
\end{align*}
We'll call the {\tt RAT} sequence $(F_k:\ k\ge 1)$ {\it compactly generated} if the collection $\{F_k:\ k\ge 1\}$ is compact.
 \proclaim{8.2 Compactness lemma}\ \ 
 Suppose that $2\a=[n_1,n_2,\dots]\in\text{\tt BAD}$ with
$$\sup_{k\ge 1}\,n_k=:M\ \text{and}\ \max\,\{J\ge 1: \ \exists\ k,\ n_{k+i}=2\ \forall\ 1\le i\le J\}=:\frak r$$
and let $(F_n:\ n\ge 1\}$ be an  $\a$-{\tt RAT}, then 
$$\frak F:=\{F_{j_1}\circ F_{j_2}\circ\dots\circ F_{j_r}:\ j_1,j_2,\dots,j_r\in\mathbb N\ \text{and}\ 1\le r\le \frak r\}$$ is compact.
\endproclaim

Before proving the compactness lemma, we need a
\proclaim{8.3 Sublemma} \ \ If $\a\in\text{\tt BAD}$, then $\exists\ \D>0$ so that
\sms{\rm (i)} $r_k:=\frac{\ell_k(1)}{\ell_k(0)}\ge\D\ \forall\ k\ge 1$,   and
\sms{\rm (ii)}  $\frac{p_k(i)}{1-p_k(i)}\in [n_{k+1}-1-i,\frac{n_{k+1}-1-i}\D]$.\endproclaim\demo{Proof of {\rm (i)}}\ \ \ \ 
\begin{align*}r_{k+1}&=\frac{\ell_{k+1}(1)}{\ell_{k+1}(0)}\\ &=\frac{(n_{k+1}-2)\ell_k(0)+\ell_k(1)}{(n_{k+1}-1)\ell_k(0)+\ell_k(1)}
\\ &=1-\frac1{n_{k+1}-1+r_k}.\end{align*}
In case $n_{k+1}\ge 3$, we have
$$r_{k+1}\ge \frac12.$$
When $n_{k+1}=2, \ r_{k+1}=v(r_k)$ where \ $v(x):=\frac{x}{1+x}.$
Now suppose that  $n_k\ge 3\ \text{and}\ n_{k+1}=n_{k+2}=\dots=n_{k+R}-2$, then for $\forall\ 1\le j\le R$,
$$r_{k+j}=v(r_{k+j-1})=v^j(r_k)\ge v^j(\tfrac12).$$
Inspection of this recursion shows that
$$\frac1{r_{k+j}}=\frac1{r_k}+j\le 2+j\ \ \forall\ 1\le j\le R.$$
Let $\frak r\in\Bbb N$ be so that $\forall\ k\ \exists\ J\le \frak r+1,\ n_{k+J}\ge 3$ and set $\D:=\frac1{2+\frak r}$.
It follows that
$$r_k\ge \D.\ \ \ \ \CheckedBox$$
\demo{proof of {\rm (ii)}}
\begin{align*}\frac{p_k(i)}{1-p_k(i)} &=\frac{(n_{k+1}-1-i)\ell_k(0)}{\ell_k(1)}\\ &=\frac{n_{k+1}-1-i}{r_k}\\ &\in
[n_{k+1}-1-i,\frac{n_{k+1}-1-i}{\D}].\ \ \ \CheckedBox
\end{align*}

\demo{Proof of the compactness lemma}
\

 We claim first that the {\tt RAT} collection $\mathcal F_0:=\{F_n:\ n\ge 1\}$ is compact. 
 \
 
 Condition (ii) is immediate from the boundedness of $\{n_k:\ k\ge 1\}$.
 To see (i), note first that \bul $q_2=1\ \text{and}\ q_N\in [\tfrac13,\tfrac12]\ \forall\ N\ge 3$, \bul $p_k(1)=0$ when $n_k=2$ and by the sublemma 
 \bul $\exists\ R>0$ so that $p_k(i),\ 1-p_k(i)\ge R\ \forall\ k\ge 1,\ i=0,1,\ (n_k,i)\ne (2,1)$.
 \
 
 Next, for each fixed $k\ge 1,\ i,j=0,1\ \text{and}\ \nu\in\mathbb Z$,
$$P([a_{i,j}(F_k)=\e]\cap [b_i(F_k)=\nu])$$ is a  polynomial of degree at most $3$ in 
$$\u z^{(k)}=(z_1^{(k)},\dots,z_9^{(k)})=(p_k(0),\ 1-p_k(0),\ q_{n_k},\ 1-q_{n_k},p_k(1),\ 1-p_k(1),\ q_{n_k-1},\ 1-q_{n_k-1},\frac1{n_k})$$
whose coefficients are non-negative integers. By sublemma 8.3, $\exists\ \e>0$ so that
$$1\le s\le 9,\ k\ge 1,\ z_s^{(k)}>0\ \implies\ z_s^{(k)}\ge \e.$$
Condition (i) follows because if
$$P([a_{i,j}(F_k)=\e]\cap [b_i(F_k)=\nu])=F(\u z)=\sum_{1\le r,s,t\le 9}N_{r,s,t}z_rz_sz_t>0$$ then
$\exists\ 1\le r,s,t\le 9$ (maybe not distinct) so that $N_{r,s,t}\ge 1$ and $z_r,\ z_s,\ z_t\ge\e$.
This shows that $\mathcal F_0$ is compact.

\

For analogous reasons, the compactness persists among concatenations of $\{F_n:\ n\ge 1\}$  of length bounded by
 $\frak r$. 

\

Let $F=F_{k_r}\circ F_{k_{r-1}}\circ\dots\circ F_{k_1}\in\frak F$, then 
$$a(F)=a(F_{k_r})a(F_{k_{r-1}})\dots a(F_{k_1})\ \text{and}\ b(F)=\sum_{\nu=1}^{r-1}a(F_{\nu+1}^r)b(F_{k_\nu})$$
where $F_{\nu+1}^r:=F_{k_r}\circ F_{k_{r-1}}\circ\dots\circ F_{k_{\nu+1}}$.

\

Condition (ii) is immediate since
$$\|b(F)\|_\infty\le \sum_{\nu=1}^{r-1}\|a(F_{\nu+1}^r)b(F_{k_\nu})\|_\infty\le \sum_{\nu=1}^{r-1}\|b(F_{k_\nu})\|_\infty\le \frak r\sup_n\|b(F_n)\|_\infty.$$

\

To see (i) we note  that for fixed $i,j=0,1\ \text{and}\ \nu\in\mathbb Z$,
$$P([a_{i,j}(F)=\e]\cap [b_i(F)=\nu])$$ is a  polynomial of degree at most  most $r$ in the variables 
$$\{P([a_{k,\ell}(F_{k_u})=\om]\cap [b_k(F_{k_u})=\mu]):\ u=1,2,\dots,r,\ \om=0,\pm 1,\ \mu\in\mathbb Z\}$$
whose coefficients are non-negative integers and so 
\begin{align*}P([a_{i,j}(F)=\e]\cap &[b_i(F)=\nu])>0\ \Lra\ \\ &P([a_{i,j}(F)=\e]\cap [b_i(F)=\nu])\ge\frac1{M_{\frak F_0}^{\frak r}}
\end{align*}
where $M_{\frak F_0}$ is as in the definition of compactness of $\mathcal F_0$.
This proves the compactness lemma.\ \ \ \Checkedbox
\subsection*{Proof of the Grouping theorem}

\

 In order to prove that  $(G_n:\ n\ge 1)$ has {\tt BMF}, it suffices by lemma 7.3 and the compactness lemma, to show that each  $G_n$ is mean contractive.
 \
 
We claim also that in order to prove adaptedness of the sequence $(G_n:\ n\ge 1)$, it suffices to  show that each {\tt RAT} $G_n$ is adapted.
\

To prove  this latter claim, suppose that each $G_n$ is individually adapted.  Compactness of $(G_n:\ n\ge 1)$ entails  $\|b(G_n)\|\le M$ which implies by discreteness that 
(possibly increasing $M$)
$\forall\ k\ge 1,\ \text{\tt supp}\,\,b_k(G_n)\subset [-M,M]$ whence 
$$\G_{G_n}\subset\frac{2\pi}{(2M)!}\mathbb Z.$$
Adaptedness of the sequence now indeed follows individual adaptedness by  compactness.

\

The rest of this proof is concerned with establishing these individual properties for the concatenations involved.
\subsection*{Concatenations of {\tt spec-RAT}s}\ \ Each $G=G_n$ in the canonical grouping is an independent concatenation of four {\tt RAT}s of form
$$H=R\circ T$$
where $R,\ T$ are independent, $R=R(H)$ is a non-trivial {\tt spec-RAT} and either $T=T(H)$ is a concatenation of finitely many independent, trivial {\tt spec-RAT}s, 
or
$T=\text{\tt Id}$.
\

\sms{\bf Mean contractivity of {\tt spec-RAT} concatenations}
\

Non-trivial, odd {\tt spec-RAT}s are mean contractive. Even {\it spec-RAT}s are not.
Suppose $G=H_1\circ H_2$ where and $R(H_i)\ \ i=1,2$ are  even. By the parity transition laws, each $R(H_i)$ is followed and preceded by 
odd (here necessarily trivial) {\tt spec-RAT}s $J$. It can be calculated that for $a=a(J\circ R)$ or $a=a(R\circ J)$, then 
$P(a_{i,0}=\e)>0$ for $i=0,1\ \text{and}\ \e=\pm 1$. By irreducibility $G$ is mean contractive.
\sms{\bf Adaptedness  of {\tt spec-RAT} concatenations}
\

If $G=H_1\circ H_2\circ H_3\circ H_4$ where $R(H_i)$ is even for some $i=2,3,4$, then  $R(H_i)$ is partially adapted  and $H_1\circ\dots\circ H_{i-1}$ is irreducible, 
whence  $H_1$, $H_1\circ\dots\circ H_{i-1}\circ R(H_2)$ is adapted and so is
 $G$.
\

The remaining case is where $G=H_1\circ\dots\circ H_4$ where $R(H_i)$ is odd $\forall\ 2\le i\le 4$.
\

To show that such $G$ is adapted, 
 we write 
$$G=H_1\circ U\circ T$$ where
$$U:=H_2\circ H_3\circ R(H_4).$$
\

Since $H_1$ is irreducible, it suffices to show that row $0$ is adapted for either $H_2\circ H_3$ or for $U$.
\

Suppose that row $0$ is not adapted for $H_2\circ H_3$, then, in particular the random variable $b_{0,0,1}(H_2\circ H_3)\equiv c$ is constant. 
\

For each $\a=\pm 1$, $P(J_\a)>0$ where
$$J_{\a}=[a_{0,0}(H_2)=1]\cap [a_{0,0}(H_3)=\a]\cap [a_{0,0}(R(H_4))=\a]\subset [a_{0,0}(U)=1]$$
and, on $J_{\a}$ we have that
                                
\begin{align*}
b_0(U)=b_0(H_2\circ H_3)+\a b_0(R(H_4))=c_1+\a.  
\end{align*}
Thus the random variable 
$P(b_{0,0,1}(U)=c_1+\a)>0$ for $\a=\pm 1$ and row $0$ is adapted for $U$.\ \ \Checkedbox

 This completes the proof of the Grouping theorem.\ \ \ \ \qed

\

\section*{Conclusion of the proof of the main result}
\

By the adaptedness preservation observation, the centering 
\par \ \  $(H_k:\ k\ge 1)$ of $(G_k:\ k\ge 1)$ is adapted (with the same adaptivity group).
\

By compactness, {\tt BMF} and the bounded centering lemma, $(H_k:\ k\ge 1)$ is uniform.
\

Thus $(H_k:\ k\ge 1)$  satisfies {\tt WRLLT}. Let $(X_n:\ n\ge 1)$, $(Y_k:\ k\ge 1)$ and $(Z_k:\ k\ge 1)$ be the {\tt ARW}s generated by $(F_n:\ n\ge 1)$, $(G_k:\ k\ge 1)$ and 
 $(H_k:\ k\ge 1)$ respectively, then 
 for each $i=0,1$ and $1\le p\le 2$,
\begin{align*}
 \int_\mathbb T|\widehat{\Phi}_{X^{(\nu_k)}(i)}(\th)|^pd\th &=\int_\mathbb T|\widehat{\Phi}_{Y^{(k)}(i)}(\th)|^pd\th\ \ \  \text{because the processes are identical}\\ &=
 \int_\mathbb T|\widehat{\Phi}_{Z^{(k)}(i)}(\th)|^pd\th\ \ \text{because}\  Z^{(k)}(i)=Y^{(k)}(i)-E(Y^{(k)}(i))\\ &
 \ \text{\Large $\asymp$}\ \frac1{\sqrt{k}}\ \ \text{by the {\tt WRLLT}.}
\end{align*}

By the  Visit Lemmas  2.4, 
$$\int_\mathbb T\Psi_{\ell_{\nu_k}(0)}(x)dx\ \gg\ \ \ell_{\nu_k}(1)\int_\mathbb T|\widehat{\Phi}_{X^{(\nu_k)}(1)}(\th)|^2d\th\ \asymp\frac{\ell_{\nu_k}(1)}{\sqrt k},$$
and 
$$\|\Psi_{\ell_{\nu_k}(1)}\|_\infty\le \ell_{\nu_k(0)}\int_\mathbb T|\widehat{\Phi}_{(X^{(\nu_k)}(0)}(\th)|d\th\ \asymp\frac{\ell_{\nu_k}(0)}{\sqrt k}.$$
Next, since $\a\in\text{\tt BAD}\ \ \text{and}\ \sup_{k\ge 1}(\nu_{k+1}-\nu_k)<\infty,\ \ \exists\ 1<m<M$ so that
$$m\le \frac{\ell_{\nu_{k+1}}(i)}{\ell_{\nu_{k}}(j)}\le M\ \ i,j=0,1$$
whence
$$\log\ell_{\nu_{k}}(j)\asymp k.$$
For $\ell_{\nu_k}(0)\le n\le \ell_{\nu_{k+1}}(0)$,
$$\int_\mathbb T\Psi_n(x)dx\ge \int_\mathbb T\Psi_{\ell_{\nu_k}(0)}(x)dx\ \gg\ \frac{\ell_{\nu_k}(1)}{\sqrt k}\gg\frac{\ell_{\nu_{k}+1}(1)}{\sqrt k}\gg
\frac{n}{\sqrt{\log n}}$$
and for $\ell_{\nu_k}(1)\le n\le \ell_{\nu_{k}+1}(1)$,
$$\|\Psi_n\|_\infty\ll   \frac{\ell_{\nu_{k+1}}(1)}{\sqrt k}\ll\frac{\ell_{\nu_k}(0)}{\sqrt k}\ll \frac{n}{\sqrt{\log n}}.\ \ \ \CheckedBox\text{(\dsrailways) (as on page \pageref{choochoo})}$$

\subsection*{Concluding Remarks}
\

\sms {\bf 1.}  It is shown in \cite{CR} that for a.e. $\a\in\Bbb T$, $T_\a$ is $\tfrac1n$-recurrent in the sense that
$\sum_{n=1}^\infty\tfrac1nf\circ T_\a^n=\infty$ a.e. $\forall\ f\in L^1_+$.  
\

It follows easily from (\dsrailways) (as on page \pageref{choochoo}) that for $\a\in\text{\tt BAD}$ and $\om_n\downarrow 0$:
\

$T_\a$ is $\om_n$-recurrent 
($\sum_{n=1}^\infty\om_nf\circ T_\a^n=\infty$ a.e. $\forall\ f\in L^1_+$) iff 
$$\sum_{n\ge 1}\frac{n(\om_n-\om_{n+1})}{\sqrt{\log n}}=\infty$$
e.g. $\om_n=\frac1{n\sqrt{\log n}}$.
\sms {\bf 2.} The subsequence condition  (\dstechnical) (as on page \pageref{techy}) is satisfied for some $\a\notin\text{\tt BAD}$.
For example, if  
$2\a=[2m_1+1,2m_2,2m_3,\dots]$ with $m_k\ge 2\ \forall\ k\ge 1$, then 
 by sublemma 8.3(i), $\ell_k(1)\le\ell_k(0)\le 2\ell_k(1)$
  and it follows as above that  (\dstechnical)  is satisfied along $\ell_{\mu_k}(0)$ ($\mu_k\uparrow\infty$) as soon as
$$\int_\mathbb T|\widehat{\Phi}_{X^{(\mu_k)}(0)}(\th)|d\th\ \ll\ \int_\mathbb T|\widehat{\Phi}_{X^{(\mu_k)}(1)}(\th)|^2d\th.$$

The parities $\e_k(0)=1\ \forall\ k\ge 2$, whence the associated {\tt RAT} collection $\{F_{j}:\ j\ge 1\}$ is compact. \
 
The  canonical subsequence in \S8 is 
$\nu_k=1+4k\ \ \ (k\ge 1)$, the grouping lemma applies and  the centering of the canonical grouping {\tt RAT} 
sequence satisfies the {\tt WRLLT}.
\

Consequently $T_\a$ satisfies (\dstechnical)  along $\ell_{\nu_k}(0)$ with $a_{\ell_{\nu_k}}\asymp\ \frac{\ell_{\nu_k}}{\sqrt k}$.
\

\sms {\bf 3.}  It is thus natural to ask if $p(\text{\dstechnical})=1$ (or even if  $p(\text{\dstechnical})>0$) where
$$p(\text{\dstechnical}):=m(\{\a\in (0,1):\ T_\a\ \text{satisfies\ (\dstechnical)\ along some subsequence}\}).$$

\end{document}